\def\boxit#1{\vbox{\hrule\hbox{\vrule\kern6pt
          \vbox{\kern6pt#1\kern6pt}\kern6pt\vrule}\hrule}}

\newcommand{\Fbar}{\overline{F}}

\newcommand{\Hbar}{\overline{H}}

\newcommand{\Gbar}{\overline{G}}
\newcommand{\Jbar}{\overline{J}}

\newcommand{\topr}{\stackrel{\mathrm{P}}{\longrightarrow}}

\newcommand{\eqdr}{\stackrel{\mathrm{D}}{=}}

\newcommand{\R}{\Bbb{R}}

\newcommand{\N}{\Bbb{N}}

\newcommand{\rmd}{{\rm d}}

\newcommand{\halmos}{\quad\hfill\mbox{$\Box$}}

\newcommand{\NN}{\mathbb{N}}

\newcommand{\dto}{\downarrow}
\newcommand{\uto}{\uparrow}

\newcommand{\be}{\begin{equation}}
\newcommand{\ee}{\end{equation}}
\newcommand{\bea}{\begin{eqnarray}}
\newcommand{\eea}{\end{eqnarray}}
\newcommand{\bean}{\begin{eqnarray*}}
\newcommand{\eean}{\end{eqnarray*}}
\newcommand{\ben}{\begin{equation*}}
\newcommand{\een}{\end{equation*}}
\newcommand{\ba}{\begin{align}}
\newcommand{\ea}{\end{align}}

\def\nexto{\kern -0.54em}

\def\Kless{K^u}
\def\Kgreater{K^c}

\documentclass[bernoulli,preprint]{imsart}

\usepackage{amsmath}
\usepackage{graphicx}
\usepackage[shortlabels]{enumitem}

\usepackage{amsfonts,pdfsync}
\usepackage{amsthm}

\usepackage{cases}
\usepackage{grffile}
\usepackage{caption}

\usepackage{diagbox}

\RequirePackage[OT1]{fontenc}
\RequirePackage{amsthm,amsmath}
\RequirePackage[numbers]{natbib}
\RequirePackage[colorlinks,citecolor=blue,urlcolor=blue]{hyperref}

\arxiv{2103.01337}

\startlocaldefs
\numberwithin{equation}{section}
\theoremstyle{plain}

\newtheorem{theorem}{Theorem}[section]

\newtheorem{corollary}{Corollary}[section]
\endlocaldefs

\begin{document}
\bibliographystyle{plain}
\begin{frontmatter}

\title{Splitting the Sample at the Largest Uncensored   Observation}
\runtitle{Splitting  at the Largest Uncensored   Observation}

\begin{aug}
\author{\fnms{Ross} \snm{Maller,}\corref{}
}
\author{\fnms{Sidney} \snm{Resnick}
}
\and\author{\fnms{Soudabeh} \snm{Shemehsavar}
}

\runauthor{Maller, Resnick \& Shemehsavar}


\address{Research
School of Finance, Actuarial Studies \& Statistics, \\
The Australian National University, Canberra, Australia\\
Ross.Maller@anu.edu.au\\
}

\address{ School of Operations Research \& Information Engineering, \\
Cornell University,  Ithaca, New  York, USA\\
sir1@cornell.edu\\
}

\address{School of Mathematics, Statistics and Computer Science, \\
University of Tehran, Tehran, Iran\\
shemehsavar@ut.ac.ir (corresponding author)\\
}
\end{aug}

%

\begin{abstract}
\noindent
We calculate finite sample and asymptotic distributions for the largest censored and uncensored survival times, and some related statistics,
from  a sample of survival data generated according to an iid  censoring model. These statistics  are important for  
assessing whether there is sufficient follow-up in the sample to be confident of the presence of immune or cured  individuals in the population.
 A key structural result obtained  is that, conditional on the value of  the largest uncensored survival time, and knowing   the number of censored  observations exceeding this time, 
the sample partitions into two independent subsamples, each subsample having
the distribution of an  iid  sample of censored survival times, of reduced size, from truncated random variables. 
This result provides valuable insight into the construction of censored survival data,
and facilitates the calculation of explicit finite sample formulae.
We illustrate for distributions of  statistics useful for testing for sufficient follow-up in a sample, and apply extreme value methods  to derive asymptotic distributions for some of those.
\end{abstract}



\end{frontmatter}

{}\null
\vskip-1cm 
%
%
%

\section{Introduction}\label{intro}
In the analysis of censored survival data, 
it was long considered an anomalous situation to observe
a  sample Kaplan-Meier estimator (KME) \cite{kaplan:meier:1958} 
 which is improper, i.e., of total mass less than 1.  This occurs when the largest survival time in the sample is  censored.
In some early treatments it was advocated to remedy this situation by redefining 
 the largest survival time in the sample to be uncensored
 (cf. Gill \cite{gill:1980}, p.35).
Later it was recognized that one or more of the largest observations being censored conveys important information concerning the possible existence  of immune or cured individuals in the population. One aim of the present paper is to focus attention on the importance of considering the location and conformation of the censored observations in the sample.

Even in a population in which no immunes are present,  it can be the
case that the largest or a number of the largest  survival times are
censored according to whatever chance censoring mechanism is operating
in the population.  The possible presence of immunes will be signaled
by an interval of constancy of the KME at its right hand end, with the
largest observation being censored. The length of that  interval  and
the  number of censored survival times larger than the largest
uncensored survival time  are  important statistics for  testing for
the presence of immunes, and for assessing whether there is sufficient
follow-up in the sample to be confident of their presence. We focus
particularly on the ``sufficient follow-up" issue in the present
paper. 

Our aim is to derive distributional results that lead to rigorous methods for
deciding if the population contains immunes and if the length of observation is
sufficient. We proceed by calculating finite sample and asymptotic joint distributions for the largest censored and uncensored survival times under an iid (independent identically distributed) 
censoring model for the data. These calculations are facilitated by a significant
splitting result which states that, conditional on knowing
both the value of the largest uncensored survival time and the number of censored observations exceeding this observation,
the sample partitions into two conditionally independent subsamples. The observations in each subsample have the distribution of an iid sample of censored
survival times, of reduced size, from appropriate truncated random variables.
See  Theorem \ref{th1} below for  a detailed description of this result.

Using this result we are able to calculate expressions for distributions of statistics, such as those mentioned, related to testing for  sufficient follow-up.
Those distributions are reported in Section \ref{s1}.
In Section \ref{asymp} we calculate large sample distributions for the statistics  under very general conditions on the tails of the censoring and survival distributions. 
Special emphasis is placed on the realistic situation where the censoring distribution has a finite right endpoint.
These results complement analysis in \cite{MR2020} made under the assumption of an infinite right endpoint.
That paper is concerned only with asymptotics whereas in the present paper we have the finite sample distributions in
 Section \ref{s1}.
  With these we can calculate moments of rvs, etc.
  
A discussion Section \ref{tfu} contains a motivating example outlining how the results can be applied in practice  to validate and  improve existing procedures.  
Proofs for the formulae in Section \ref{s1}  are in Section \ref{pfs} and  the asymptotic results are  proved in Section \ref{sspa}. 
A supplementary arXiv submission \cite{arxivsub}
contains plots illustrating some of the distributions in Section \ref{s1}. 

For further  background  we refer to the book by Maller and Zhou \cite{MZ1996} which gives many practical examples
from medicine, criminology and various other fields of this kind of data and its  analysis.
See also  the review articles by  
Othus,  Barlogie,  LeBlanc and Crowley \cite{oblc:2012}, 
Peng  and Taylor \cite{pt:2014}, 
Taweab and  Ibrahim \cite{ti:2014},
Amico and  Van Keilegom \cite{avk2018},
and  a recent paper by 
 Escobar-Bach, Maller, Van Keilegom and Zhao \cite{emvz:2020}.


\newpage

\section{Splitting the Sample: Finite Sample Results}\label{s1}

\subsection{Setting up: the iid  censoring model}\label{s2}
\vskip-.3cm
We assume a general independent censoring model with right censoring and adopt the following formalism. 
Assume $F(x)$ and $ G(x)$ are  two proper, continuous,  cumulative distribution functions (cdfs) 
on $[0,\infty)$.
Let  $p\in(0,1]$ be  a parameter. Define the iid 
sequence of  triples $\{(L_i, U_i,
B_i),i\geq 1\}$, where $L_i$ has distribution
$F$, $U_i$ has distribution $G$ and $B_i$ is Bernoulli with
$$P(B_i=1)=1-P(B_i=0)=p.$$
For each $i$ the  components $L_i$, $U_i,$ and $B_i$ are assumed  mutually independent.
The $U_i$ are censoring variables and the $L_i$ represent lifetimes
given that they are finite.
Define
$$
T_i^*=\infty \cdot (1-B_i) + L_iB_i
=\begin{cases}
  \infty,&  \text{ if }B_i=0,\\
  L_i, & \text{ if } B_i=1.\end{cases}
$$
 Then $T_i^* $ has distribution $F^*$ given as 
\be\label{7}
P(T_i^*\leq x)=F^*(x):=pF(x),\, 0\leq x <\infty, \ \text{and}
\ P(T_i^*=\infty)=1-p.
\ee
We allow $T_i^*$ to be infinite to cater for cured or immune individuals  in the population. 

The auxiliary unobserved Bernoulli random variable (rv)  $B_i$ indicates whether or not individual $i$ is immune;
 $B_i=1$ (resp., $0$)  indicates that  $i$ is  susceptible to death or failure (resp., immune).
 When $p=1$, equivalently $B_i\equiv 1$, all individuals are susceptible; this is the situation in ``ordinary" survival analysis (when  all individuals in the population are susceptible to death or failure).
 We do not observe the $B_i$ so we do not know whether an individual is susceptible or immune. Overall, our model is commonly  known as the ``mixture cure model".
 
Immune individuals do not fail so  their lifetimes are  formally taken  to be infinite.
Thus we set
$T_i^*=\infty$ when $B_i=0$, 
with $P(T_i^*=\infty)=1-p$ 
being the probability individual $i$ is immune. 
 The cdf $F$ can be interpreted as the distribution of the susceptibles' lifetimes.

As is usual in survival data, and essential in our study of long-term survivors, potential lifetimes are censored at a limit of follow-up represented for individual $i$ by the random variable $U_i$.
The observed lifetime for individual $i$ is thus of the form
  \be\label{olt}
  T_i:= U_i\wedge T_i^*.
\ee 
 Immune individuals' lifetimes are, thus, always censored.
The distribution of the $T_i$ is $H(x):= P(T_i^*\wedge  U_i\le x)$.
We also observe 
indicators of whether censoring is present or absent:
$$C_i=1_{\{T^*_i\leq U_i\}}.$$
Tail or survivor functions of  distributions $F$, $G$, $F^*$ and $H$  on $[0,\infty)$ are denoted by
$\Fbar=1-F$, $\Gbar=1-G$,  $\Fbar^*=1-F^*$ and $\Hbar (x)=
1-H(x)$.
We note that 
\begin{equation}
\Hbar (x)=P(T_i^*\wedge  U_i>x)=\Fbar^*(x)\Gbar(x)
=\big(1-p+p\Fbar(x)\big)\Gbar(x), \label{13} 
\end{equation}
and 
\ben
\Fbar^*(t):=P(T^*>t) = 1-p+p\Fbar(t), \ t\ge 0.
\een

To connect this formalism with the observed survival data, we think of 
the $T_i^*$ as representing the times of occurrence of an event under study, such as  the death of a person, the onset of a disease, the recurrence of a disease, the arrest of a person charged with a crime, the re-arrest of an individual released from prison, etc. 
We allow the possibility that only a proportion $p\in(0,1]$ of individuals in the population are susceptible to death or failure, and the remaining $1-p$ are ``immune" or ``cured".

\subsection{Splitting the sample: main structural result}\label{s1a}
In addition to the notation in the previous subsection, let $M(n): =\max_{1\leq i\leq n}T_i$ be the largest observed survival time and let 
$M_u(n)$  be the largest observed {\it uncensored} survival time.   Since all $T_i>0$ with probability 1 we can define 
$M_u(n)$ in terms of the $T_i$ and $C_i$  by 
$M_u(n) =\max_{1\le i\le n } C_iT_i$.

To state the splitting result, we adopt the convention that for a
non-negative random variable $X$ and Borel set $B\subset
\R_+$
 with $P(B)>0$,
 $(X|X\in B) $ is a random variable with distribution
\begin{equation}\label{e: conditDistn}
  P(X\in A|X\in B)=P(X\in A\cap B)/P(X\in B), 
  \quad A\subset \R_+,\,A
\text{ Borel}.\end{equation}
Our splitting theorem says that
the sample $S_n:=\{T_i,1\le i\le n\}$  partitions into
\begin{align}
  S_n:=\{T_i,1\le i\le n\}=&\{M_u(n) \}\cup  \{T_i:i\leq n\; \&\;
                           T_i<M_u(n)\}\nonumber\\
&  \qquad \qquad
  \cup 
                           \{T_i:i\leq n \;\& \; T_i>M_u(n)\}\nonumber\\
  =: &\{M_u(n)\} \cup S_n^< \cup S_n^> . \label{e:partition}
\end{align}
Conditional on knowing that $ M_u(n)=t>0$ and that $r$  censored  observations exceed $ M_u(n)$, $0\le r\le n-1$, 
$S_n^<$ consists of $n-r-1$ iid variables with distribution that of $(T_i|T_i<t)$,
     and $S_n^>$ consists of $r$ iid variables with tail function
     \begin{equation}\label{e:censort_tail}
P(T_i^{>,c}(t) >x):=    \frac{\int_x^\infty \Fbar^*(s) G(ds)}{\int_t^\infty \Fbar^*(s)
       G(ds)},\quad x\geq t,\end{equation}
     which is the conditional distribution tail of a censored observation
     given that it is bigger than $t$. 
     (See Eq. (2.13)  in \cite{MR2020}.)
     Furthermore, $S_n^<$ and $S_n^>$ are conditionally independent. Note that observed lifetimes less
     than $M_u(n)$ may be either censored or uncensored but observed
     lifetimes greater than $M_u(n)$ are necessarily censored.

For further  precision in stating the splitting result, we need
notation for the numbers of  censored  observations
smaller or greater than  $M_u(n)$,  
the largest uncensored survival time in the sample.
On $\{M_u(n)>0\}$, let
\be\label{Nd1}
N_c^>(M_u(n)):= \{{\rm number\ of\ censored\ observations\ exceeding}\  M_u(n)\}.
\ee
By convention we set
\ben 
\{N_c^>(M_u(n))=0\}
=\{M_u(n)= M(n)\}
=
\{{\rm largest\ observation\ uncensored} \}
\een
and
\ben 
\{N_c^>(M_u(n))=n\}= 
\{{\rm all}\ n\ {\rm observations\ censored} \} = \{C_1=\cdots=C_n=0\}.
\een
On $\{N_c^>(M_u(n))=n\}$ we set $M_u(n)=0$.
Next, let
\be\label{Nu4}
N_u(n)
:= 
\{{\rm total\ number\ of\ uncensored\ observations\ in\ the\ sample} \},
\ee
and when $N_u(n)>1$, define 
\begin{align}\label{Nd4}
&
N_u^<(M_u(n))\cr
&:= 
\{{\rm number\ of\ uncensored\ observations\ strictly \ less\ than}\  M_u(n)\}
\end{align}
and 
\be\label{Nd5}
N_c^<(M_u(n)):= \{{\rm number\ of\ censored\ observations\  less\ than}\
 M_u(n)\}.
\ee
On $\{N_u(n)=1\}$, set $N_u^<(M_u(n))=N_c^<(M_u(n))=0$.
When $N_u(n)=0$, we do not define
$N_u^<(M_u(n))$ or $N_c^<(M_u(n))$ (and of course there's no need for a notation like $N_u^>(M_u(n))$
since such a number would always be 0.)
Lastly, let
\be\label{Nd6}
N_c(n):
= 
\{{\rm total\ number\ of\ censored\ observations\ in\ the\ sample} \}.
\ee

We also use the notation $\NN_n:=\{1,2,\ldots,n \}$, $n=1,2,\ldots,$.With these definitions and conventions,
on $\{N_u(n)\ge 1\}$ the 
$N_u^<(M_u(n))$, $N_c^<(M_u(n))$ and $N_c^>(M_u(n))$
take values in
$\N_{n-1}\cup\{0\}$,
satisfying
$N_u^<(M_u(n))+N_c^<(M_u(n))+ N_c^>(M_u(n))=n-1$
and $N_c(n)= N_c^>(M_u(n))+N_c^<(M_u(n))$, 
and we have
\ben
\{N_u(n)=0\}= 
\{{\rm all}\ n\ {\rm observations\ censored} \} =
\{M_u(n)=0\}.
\een
The precise statement of the main splitting result is contained in the next theorem and proved in Section \ref{pfs}. Theorem \ref{th1} also contains the  distribution in \eqref{j0} needed for later calculations.
Let  $\tau_{F}= \sup\{t>0:F(t)<1\} $ be the right endpoint of the support of the cdf $F$, and similarly define $\tau_G$, $\tau_H$ and $\tau_{F^*}$.   

\begin{theorem}\label{th1}
  [Splitting the Sample at $M_u(n)$] 
  For a sample of size $n$,

(i)  The joint distribution of 
$\bigl(M_u(n),N_c^>(M_u(n))\bigr) $ is, for $0\leq r\leq n-1$, $0\le t\le \tau_H$,
\bea\label{j0}
&&
P\big(0\le  M_u(n)\le t,\,  N_c^>(M_u(n))=r\big)\cr
&&=
n{n-1\choose r}
\int_{y=0}^t
\Big( \int_{z=y}^{\tau_H} \Fbar^* (z)\rmd G(z) \Big)^{r}
H^{n-r-1}(y)   \Gbar(y)\rmd F^*(y).\ \
\eea

(ii) Let $A_r$ and $B_{n-r-1}$ be Borel subsets of $\R_+^r$ and $\R_+^{n-r-1}$ respectively. 
Conditional on $N_c^>(M_u(n))=r$, the sets $S_n^<$ and
$S_n^>$ defined in \eqref{e:partition} have cardinality $n-r-1$ and
$r$. Then,  for $0\le t\le \tau_H$,
\begin{align}\label{w10b}
 & P\Big((T_i:T_i \in S_n^<)\in B_{n-r-1}, \,
  (T_i:T_i \in S_n^>)\in   M_u(n)+A_r           \cr
  &\qquad \Big|
    M_u(n)=t, N_c^>(M_u(n))=r\Big)\cr
  =&
  P\big((T_j|T_j<t)_{j=1,\dots,n-r-1}\in B_{n-r-1}\big) P\big((T_l^{>,c}(t))_{l=1,\dots,r}
 \in t+A_r\big),\ \ 
\end{align}
where $(T_j|T_j<t)$ is defined using   \eqref{e: conditDistn}
and $T_l^{>,c}(t)$  has the distribution tail in \eqref{e:censort_tail}.\end{theorem}
 
 \noindent{\bf Remarks.}\
(i) \ 
The case $r=0$ in Theorem \ref{th1} means no censored observation exceeds $M_u(n)$,
so $M_u(n)=M(n)$,
and the second component in the statement of the theorem is empty.

When $r=n-1$,  the first component in the statement of the theorem is empty
and all observations are censored except for the smallest which is  $M_u(n)$.

We can also include the case $r=n$ in which all observations are censored in 
 Theorem \ref{th1} if we interpret the distribution as the conditional distribution given 
$M_u(n)=0$, and again the first component in the statement of the theorem as being empty. For the specific formula, see \eqref{c0a} in Section \ref{pfs}.

(ii)\ The splitting formula \eqref{w10b} remains true if conditioning is done on
$(M_u(n), M(n), N_c^>(n))$ 
rather than on $(M_u(n),N_c^>(n))$.
This is because the extra information in $M(n)$ beyond that in $M_u(n)$ involves only  the (censored) observations greater than $M_u(n)$. 
For proof see Corollary \ref{c1}. 
\halmos

\subsection{Finite sample distributions of the maximal times}\label{ssT}

 Theorem \ref{th1} gives a way to think about censored survival data, especially in the presence of a possible immune component, and provides
a route to calculating descriptor 
distributions such as the joint finite sample distribution
of  $M(n)$ and $M_u(n) $.
 Recall that  $\tau_{H}$ is the right endpoint of the support of $H$.


\begin{theorem}\label{th4a}[Distributions of  $M_u(n)$ and $M(n) $]

(i)\ The joint distribution of  $M_u(n)$ and $M(n) $ is given by
\bea\label{JD}
&&
P\big(0\le  M_u(n)\le t,\,  0\le  M(n)\le x \big)\cr
&&=
\begin{cases}
\big(\int_{z=0}^{x}  \Fbar^* (z)\rmd G(z)\big)^n,\ t=0,\,   0\le x\le \tau_H;  \\
H^n(x),
\ 0\le t\le \tau_H;\, 0\le x\le t; \\
\big(\int_ {z=t}^{x}   \Fbar^* (z)\rmd G(z)+H(t)\big)^n,
\  0\le t<x \le \tau_H.
\end{cases}
\eea

(ii)\ The distribution of $M_u(n) $ is given by
\be\label{4f}
P\big(M_u(n) \le  t\big) =J^n(t),\ t\ge 0,
\ee
where $J(t)$ is the distribution of an uncensored lifetime:
\be\label{4g}
  J(t)= 
1-\int_{z=t}^{ \tau_H} \Gbar (z)\rmd F^*(z)=\int_ {z=t}^{\tau_H}   \Fbar^* (z)\rmd G(z)+H(t), \  0\le t\le \tau_H.
\ee
%
%
\end{theorem}

\medskip\noindent{\bf Remarks.}\
(i)\ 
There is no probability mass outside the region $[0,\tau_H]\times [0,\tau_H]$ so the distribution in  \eqref{JD} equals 1 for 
$ t>\tau_H$, $x>\tau_H$.
Likewise the distribution in  \eqref{4g} equals 1 for 
$ t>\tau_H$.

Note also that Lines 2 and 3 
on the RHS of \eqref{JD} include the value for $t=0$;
there is  mass on the interval  $\{t=0\}\times [0\le x\le \tau_H]$, as given by the first line on the RHS of \eqref{JD}. 
Illustrative plots of  the distributions of  $M(n)$ and $M_u(n)$
are in the Supplementary Material to the paper.

(ii)\ 
$M_u(n) $ has the distribution of the maximum of $n$ iid  copies of a rv with distribution $J$ on $[0,\infty)$.
The distribution  has mass  $\big(\int_{z=0}^{\tau_H}  \Fbar^* (z)\rmd G(z)\big)^n$  at 0  corresponding to all observations being censored.
(It may seem pedantic to include these degenerate cases but they are important for checking that distributions are proper (have total mass 1).)

The right extreme $\tau_J$ of the distribution $J$ may be strictly  less than $\tau_G$; in fact, we have $\tau_J=\tau_{F}\wedge \tau_G$, as is derived in the proof of Theorem \ref{th4a}.
No uncensored observation, including the sample maximum of the uncensored observations, can exceed the smaller of $\tau_{F}$ and $\tau_G$.
Note that, in general, $\tau_J\ne \tau_H=\tau_{F^*}\wedge \tau_G$. 
We always have $H(\tau_H)=1$, $G(\tau_G)=1$  and $F(\tau_{F})=1$;
when $p=1$, so that $F^*\equiv F$, then $F^*$ has total mass 1 and $\tau_{F^*}=\tau_{F}$;
when $p<1$ we have $\tau_{F^*}=\infty$, and  $\tau_{F}\le \tau_{F^*}$,
with the possibility that  $\tau_{F}< \tau_{F^*}$.

(iii)\  
$M(n)$ has the distribution of the maximum of $n$ iid  copies of a rv with distribution $H$ on $[0,\tau_H]$; namely, 
\be\label{4h}
P(M(n)\le x)=H^n(x),\ 0\le x\le \tau_H.
\ee
%

Also important for statistical purposes are the length of the time interval between  the largest uncensored survival time  and the largest survival time, and the ratio of those times. For them
we have the following distributions. 
 
\begin{theorem}\label{th4}[Distributions of $M(n) - M_u(n) $
and  $M(n)/M_u(n) $]

We have for $0\le u\le \tau_H$
\bea\label{12a}
&&
P\big( M(n)  - M_u(n)  \le u \big)= \cr
&&
n\int_{t=0}^{\tau_H}
 \Big(\int_{z=t}^{\min(t+u, \tau_H)} \Fbar^*(z) \rmd G(z) +H(t)\Big)^{n-1} \Gbar(t) \rmd F^*(t)\cr
 &&\hskip5cm 
+
 \Big(\int_{z=0}^{u} \Fbar^*(z) \rmd G(z)\Big)^n,
\eea
with $P\big( M(n)  - M_u(n)  \le u \big)= 1$ for $u>\tau_H$.
We have for $v\ge 1$
\bea\label{12b}
&&
P\big( M(n) \le v M_u(n) | M_u(n) >0 \big)\cr
&&\hskip0.5cm 
=
\frac{\int_{t=0}^{\tau_H}
 \Big(\int_{z=t}^{\min(tv, \tau_H)} \Fbar^*(z) \rmd G(z) +H(t)\Big)^{n-1} \Gbar(t) \rmd F^*(t)}
{\int_{t=0}^{\tau_H}
 \Big(\int_{z=t}^{\tau_H} \Fbar^*(z) \rmd G(z) +H(t)\Big)^{n-1} \Gbar(t) \rmd F^*(t)},
\eea
with $P\big( M(n)  \le v M_u(n) | M_u(n) >0  \big)= 0$ for $0\le v<1$.
Both \eqref{12a} and \eqref{12b} remain true for $0<\tau_H\le \infty$. 
\end{theorem}

\medskip\noindent{\bf Remarks.}\
(i) \ 
Setting $u=0$ in \eqref{12a}  we see that the distribution of the difference  $M(n) - M_u(n) $ has mass at 0 of
\be\label{12c}
P\big( M(n)  - M_u(n)  =0 \big)
=
P\big( M(n)  = M_u(n)   \big)
=
n\int_{t=0}^{\tau_H}H^{n-1}(t) \Gbar(t) \rmd F^*(t), 
\ee
while setting $u=\tau_H$  in \eqref{12a} and observing that 
\be\label{F^*GH}
\rmd H(t)= \Fbar^*(t)\rmd G(t)+\Gbar(t)\rmd F^*(t),
\ee
we can check that the total mass is 1 by doing the integration 
\bea\label{12d}
&&
n\int_{t=0}^{\tau_H}
 \Big(\int_{z=t}^{\tau_H} \Fbar^*(z) \rmd G(z) +H(t)\Big)^{n-1} \Gbar(t) \rmd F^*(t)\cr
&&
=
 \Big(\int_{z=t}^{\tau_H} \Fbar^*(z) \rmd G(z) +H(t)\Big)^{n} \Big|_{t=0}^{\tau_H}
=
H^n(\tau_H) -   \Big(\int_{z=0}^{\tau_H} \Fbar^*(z) \rmd G(z)\Big)^n
 \cr
&&=
1-   \Big(\int_{z=0}^{\tau_H} \Fbar^*(z) \rmd G(z)\Big)^n.
\eea
Taking $u=\tau_H$ in the second term on the RHS of \eqref{12a} 
and adding this to the RHS of \eqref{12d}  we get 1.

(ii)\ 
The denominator in  \eqref{12b} multiplied by $n$ 
is the expression on the LHS of \eqref{12d} and equal to  $P(M_u(n)>0)$ as can be seen from \eqref{4g}.



\subsection{Distributions of the Numbers}\label{ssN}
The results in Subsection \ref{ssT} are obtained in Section \ref{pfs} as special cases of the formulae for the joint distributions of  $M(n)$, $M_u(n) $ and $N_c^>(M_u(n))$ which we derive there.
That analysis can be expanded to  obtain more generally the 
joint distribution of  $M(n)$, $M_u(n) $, $ N_c^>(M_u(n))$ and $N_c^<(M_u(n))$
(and then $ N_u^<(M_u(n))=n-1 -N_c^>(M_u(n))-N_c^<(M_u(n))$).
This allows derivation of the joint distribution of 
$N_c^>(M_u(n))$, $N_c^<(M_u(n))$ and $N_u^<(M_u(n))$, variables which are also useful in addressing questions of sufficient follow-up. 

We omit the details of this more general analysis  here,    
but give a main result concerning
the vector $(N_c^>(M_u(n)), N_c^<(M_u(n)), N_u^<(M_u(n)))$.
This vector is not as might be thought at first multinomially distributed, but it is,  conditional on the value of $M_u(n)$.
This again illustrates the simplicity of exposition gained by conditioning on $M_u(n)$. For this result, we need some more notation.
Define the functions
\bea  \label{pdefs}
p_c^>(t)
&=&
 \frac{\int_{y=t}^{\tau_H} \Fbar^* (y)\rmd G(y)}
{ \int_{y=t}^{\tau_H} \Fbar^* (y)\rmd G(y)+H(t)}, \cr
p_c^<(t)
&=&
 \frac{\int_{y=0}^{t} \Fbar^* (y)\rmd G(y)}
{ \int_{y=t}^{\tau_H} \Fbar^* (y)\rmd G(y)+H(t)},\cr
p_u^<(t)
&=&
\frac{ \int_{y=0}^{t} \Gbar (y)\rmd F^*(y) }
{ \int_{y=t}^{\tau_H} \Fbar^* (y)\rmd G(y)+H(t)},
\eea
which are non-negative and (using \eqref{F^*GH}) add to 1 for each $t\in (0,\tau_H)$. 

\begin{theorem}\label{th2}[Distributions of Numbers]

(i)\ 
We have for  $t>0$,  $0\le r,s,k\le n-1$, 
$r+s+k=n-1$,
the multinomial probability 
\bea\label{w15a}
&&
P\big(N_c^>(M_u(n))=r,\,  N_c^<(M_u(n))=s,\,  N_u^<(M_u(n))=k
 \big| M_u(n)=t  \big) \cr
&&\cr
&&\hskip2cm 
=  \frac{(n-1)!   }{r!\, s!\, k!} \times
 (p_c^>(t))^r  (p_c^<(t))^s   (p_u^<(t))^{k}.
\eea

(ii)\ 
Consequently,   
 conditional on $ M_u(n)=t$,
 the marginal rvs
  $N_c^>(M_u(n))$,
 $N_c^<(M_u(n))$ and $N_u^<(M_u(n))$
 are binomial with $n-1$ as the number of trials and success probabilities $p_c^>(t)$, $p_c^<(t)$ and $p_u^<(t)$ respectively.
 
(iii)\ Conditional on $ M_u(n)=t$, the number of censored observations
 $N_c(n)= N_c^>(M_u(n))+N_c^<(M_u(n))$
is Binomial $(n-1, p_c(t))$, 
where $p_c(t)=p_c^<(t)+p_c^>(t)$.  
 
 (iv)\  Conditional  on $N_c(n)=\ell$ and $ M_u(n)=t$, the number
 $ N_c^>(M_u(n))$ is Binomial $(\ell, p_c^+(t))$, 
 where
\be\label{pc+}
p_c^+(t):= 
 \frac{\int_{y=t}^{\tau_H} \Fbar^* (y)\rmd G(y)}
{ \int_{y=0}^{\tau_H} \Fbar^* (y)\rmd G(y)}.
\ee
\end{theorem}

\noindent{\bf Remarks.}\ Note that we keep $t>0$ in Theorem \ref{th2}, so conditioning on $M_u(n)=t$ as we do implies 
$M_u(n)>0$, thus $N_u(n)\ge 1$,  and there is at least one uncensored observation.
Thus $N_u^<(M_u(n))+N_c^<(M_u(n))+N_c^>(M_u(n))=n-1$.

\section{Asymptotic Results}\label{asymp}
In practice, samples of survival data can be large enough that asymptotic methods are appropriate. Equations 
 \eqref{4f} and \eqref{4h} suggest the use of extreme value methods to find limiting  distributions of  $M(n)$ and $ M_u(n)$.
For applications we are particularly interested in the cases when $G$ and/or $F$ have finite right endpoints.
In the theorem that follows we assume Type III Weibull extreme value domain of attraction conditions on $F$ and $G$, where the extreme value shape parameter is negative
(as well as a  tail balancing condition in Case 2 of the theorem). We refer to \cite{deHaan:Ferriera} and \cite{resnickbook:2008}
for general extreme value theory.

\begin{theorem}\label{FGu}[Asymptotic distribution of  
$(M(n), M_u(n))$]
Recall we  assume that $F$ and $G$ are continuous distributions.
We have the following limiting distributions in cases of interest. \newline
\noindent {\bf Case 1:}\ Assume  $\tau_{F}< \tau_G<\infty$ and $0<p<1$,
so that   $\tau_J=\tau_{F}< \tau_H= \tau_G<\tau_{F^*}=\infty$.
Suppose  in addition that, as $z\dto 0$, $0<z<\tau_G$, 
\be\label{F^*GD}
\Gbar(\tau_G-z)=a_G (1+o(1)) z^\gamma L_G(z) 
\ {\rm and }\
  \Fbar(\tau_{F}-z)= a_{F}(1+o(1)) z^\beta  L_{F}(z), 
\ee
where $a_G, a_{F}, \gamma, \beta$ are positive constants and $L_G(z)$ and $L_{F}(z)$ are slowly varying as $z\dto 0$. 
Then an asymptotic independence property holds for the random variables
$M(n)$ and $M_u(n)$, namely, for $u,v\ge 0$,
\begin{align}\label{as1}
\lim_{n\to\infty}P\big(a_n(\tau_G- M(n))\leq u,\, &b_n(\tau_{F}- M_u(n)) \leq   v\big)   \nonumber\\
=&
\big( 1-  e^{- (1-p) u^\gamma}\big)
\big(1  - e^{- p\Gbar(\tau_{F})v^\beta}\big),
\end{align}
for some deterministic norming sequences $a_n\to\infty$ and $b_n\to\infty$, as $n\to\infty$.

\noindent {\bf Case 2:}\ Assume  $\tau_{F}<\tau_G<\infty$ and $p=1$,
so that $F^*\equiv F$ and  $\tau_J=\tau_{F}= \tau_{F^*}=\tau_H<\tau_G<\infty$.
Suppose  in addition that, as $z\dto 0$,  $0<z<\tau_G$, 
\be\label{F^*GD2}
\Gbar(\tau_G-z)=a (1+o(1)) z^\beta L(z) 
\ {\rm and }\
  \Fbar(\tau_{F}-z)= a(1+o(1)) z^\beta  L(z), 
\ee
where $a$ and $\beta$ are positive constants and $L(z)$ is slowly varying as $z\dto 0$. 
Then   there exists a  deterministic sequence $a_n\to\infty$ as $n\to\infty$ such that, for $ u,v\ge 0$,
\be\label{as31}
\lim_{n\to\infty}
P\big( a_n(\tau_{F}- M(n))\le u,\,  a_n(\tau_{F}- M_u(n))\le v \big)
=1-  e^{- \Gbar(\tau_{F}) u^\beta}. 
\ee

\noindent {\bf Case 3:}\ Assume  $\tau_G<\tau_{F}<\infty$ and $0<p<1$,
so that   $\tau_J= \tau_H= \tau_G<\tau_{F}<\tau_{F^*}=\infty$,
and assume the first relation in
\eqref{F^*GD} holds.
Suppose  in addition that,  in a neighbourhood of $\tau_G$, $F$  has a density $f$ which is positive and continuous at $\tau_G$. 
Then   there exist  deterministic sequences $a_n\to\infty$ and $b_n\to\infty$ as $n\to\infty$ such that
$a_n(\tau_G- M(n))$ and $b_n(\tau_G- M_u(n))$ are asymptotically independently  distributed with 
marginal   cdfs,  respectively, 
\bea\label{as2}
&&
 1-  e^{-(1-pF(\tau_G))u^\gamma}
\ {\rm and}\ 
1  - e^{-p f(\tau_G) v^{1+\gamma}/(1+\gamma)}, \ u,v\ge 0.
\eea
Further, this  result remains true under the same assumptions
when $p=1$, and/or  when $\tau_{F}=\infty$.   
%

\noindent 
In any of Cases 1--3 we have  $M(n)\topr \tau_H$ and $ M_u(n)\topr \tau_J$ as $n\to\infty$.    
\end{theorem}

\noindent{\bf Remarks.}\
(i)\ 
Under the assumptions of  Theorem \ref{FGu}, 
 $a_n(\tau_G- M(n))$ and $b_n (\tau_{F}- M_u(n))$ (in Case 1) 
 or  $b_n (\tau_G- M_u(n))$ (in Case 3) 
 are asymptotically independent with Weibull distributions, or, as a special case, exponential. The required choices of $a_n$ and $b_n$ are specified in the proof of the theorem.

(ii)\ 
A common  assumption is of an exponential distribution for  lifetime survival: 
$F(t)= 1-e^{-\lambda t}$, $\lambda>0$, $t\ge 0$, and the uniform distribution for censoring,  $G=U[A,B]$.
This situation,  or a close approximation to it, is often the case in practice. 
See for example  Goldman \cite{gold84}, \cite{gold91}, who assumes a scenario of patients being accrued to a trial at random times for a fixed period. The survival distribution is assumed exponential and the censoring is uniform over a known period. 
These distributions for $F$ and $G$ constitute  very good baseline reference distributions for assessing the practicality of  theoretical results. 

(iii)\
When $G=U[0,\tau_G]$, $\tau_G>0$,  we have
 $\Gbar(\tau_G-z) = (z/\tau_G){\bf 1}_{\{0\le z\le \tau_G\}}$.
 Thus $G$ satisfies \eqref{F^*GD} with $a_G=1/\tau_G$, $\gamma=1$ and $L_G\equiv 1$, while $\tau_{F}=\infty$ when $F$ is exponential$(\lambda)$. 
 Case 3  of Theorem \ref{FGu} applies.

(iv)\  For a complementary  asymptotic analysis of $M_u(n)$ and $M(n)$ when $F$ and $G$ have infinite right endpoints and are in the domain of attraction of the Gumbel distribution, see  \cite{MR2020}. 
We note that the Gumbel domain also includes distributions with finite end point.

(v)\ In Case 3 of Theorem \ref{FGu}, the requirement that
$F$  have a positive density in a neighbourhood of $\tau_G$ can be replaced with the less restrictive  assumption that the quantity 
$F(\tau_G-1/x, \tau_G]$, $x>1/\tau_G$
(the mass assigned by $F$ to the interval $(\tau_G-1/x, \tau_G]$) is of the form $c_\delta x^{-\delta}L(1/x)$, where
$\delta>0$, $c_\delta>0$ and $L(1/x)$ is slowly varying as $x\to\infty$.
The factor $f(\tau_G)/(1+\gamma)$ is then replaced with
$c_\delta \delta/(\delta +\gamma)$ in \eqref{as2}.
The version in \eqref{as2} is recovered under the assumptions of the theorem when $c_\delta$ is replaced by $f(\tau_G)$ and $\beta=1$. We omit the details.

\subsection{Exact vs Asymptotic}\label{s43}
Here we give a small illustration of how the results can be used. 
We concentrate on Case 3 of Theorem \ref{FGu}, the case of 
insufficient follow-up. 
Table \ref{eqtable1} has the 95\% quantiles of 
the distribution of  $M (n) - M_u(n)$    in \eqref{12a}
 assuming  uniform censoring and a unit  exponential survival distribution truncated at $\tau_F= 4.61$ 
 (a unit  exponential  with   99\% of its mass below $\tau_F$).
Sample sizes of $n=50, 100, 500, 5000, 20000$, are listed, with  $n=\infty$ denoting the corresponding quantiles from  the asymptotic distribution in \eqref{as2}.
Values of $\tau_G$ range from 1 (very heavy censoring, insufficient follow-up) to $\tau_G=4$ (lighter censoring, but still insufficient follow-up).
Susceptible proportion is $p=0.7$. 

\begin{table} [hb]  
\begin{center}
\begin{tabular}{c|cccccc}
 \hline
\diagbox{$\tau_G$}{$n$}&50&100&500 &5000 &20000&$\infty$  \\
  \hline
    \hline
$1$&0.528&0.401&0.198&0.067&0.034&0\\
$2$&1.004&0.796&0.425&0.151&0.079&0\\
$3$&1.596&1.317&0.760&0.292&0.164&0\\
$4$&2.293&1.954&1.224&0.438&0.319&0\\
  \hline
  $1$&3.73&4.01&4.43&4.74&4.81&5.77\\
  $2$&5.02&5.62&6.72&7.55&7.90&9.51\\
$3$&6.52&7.59&9.81&11.92&13.39&15.67\\
$4$&8.11&9.77&13.68&15.49&22.56&25.84\\
  \hline
    \hline
\end{tabular}
\caption{95\% quantiles for $M(n)-M_u(n)$.\\
Upper panel, unscaled, lower panel, scaled. }\label{eqtable1}
\end{center}
\end{table}

The upper panel in Table \ref{eqtable1} has the unscaled 95\% quantiles; in the lower panel $M (n) - M_u(n)$ 
is scaled by the factor giving the asymptotic distribution in  Theorem \ref{FGu} (the second distribution in  \eqref{as2}).
The approach to  the asymptotic distribution is rather slow
but faster when censoring is heavy (small $\tau_G$).
In a sample situation we could use estimated distributions in place of those assumed.


\section{Discussion: Testing for Sufficient follow-up}\label{tfu}
Maller and Zhou \cite{MZ1996} coin the phrase  {\it sufficient follow-up}  and, based on properties of the KME, specify it to be present in a population when $\tau_{F}\le \tau_G$. (A rationale for this is in Sections 2.2 and 2.3 of   \cite{MZ1996}). To perform a test of this, we assume the contrapositive hypothesis, $H_0: \tau_G<\tau_{F}$. This implies that $\tau_G<\infty$, and for this reason we
have emphasised  the necessity to include distributions with finite right endpoints in our analyses. 
 Cases 3 and 4 of Theorem \ref{FGu} are especially relevant in this context.
Notwithstanding this, distributions with infinite right endpoints (such as, for example, the exponential, Weibull or lognormal) are commonly used in practice to model survival times (and, possibly the censoring times). There may indeed be situations in which an  infinite right endpoint is not unrealistic (certain individuals can be strictly immune to a disease, in that they can never catch the disease);
 in other situations the assumption is a close enough approximation for practical purposes.
For a comprehensive theory  we need to thoroughly explore both possibilities, as in
\cite{emvz:2020} and  \cite{MR2020}.

Statistics for  testing $H_0$
can be constructed from 
 the  number of censored survival times larger than the largest uncensored survival time  
and/or the length of the interval of constancy of the KME at its right hand end.
Consequently we have concentrated on finding the distributions of these and similar  quantities as an aid to the development of   rigorous statistical methods.
Apart from  providing foundational results for this purpose, Theorem \ref{th1} gives  strong intuitive insight into the structure of a censored sample.

The importance of developing a reliable  test for sufficient follow-up is underscored by
a recent  study of Liu et al. \cite{Liu2018},  
 in which such  testing is done on a very extensive scale.
Those researchers processed  follow-up data files for 11,160 patients
across 33 cancer types, 
calculating median follow-up times as well as median times
to event (or censorship) based on the observed times for 
four endpoints (overall survival, disease-specific survival, disease-free interval, or progression-free interval).
They classified all $33 \times 4$ resulting  KMEs  as having sufficient or insufficient follow-up (or noted cases in which  tests were inconclusive)  in order 
to give  endpoint usage recommendations for each cancer type.
They stress: {\it 
For each endpoint, it is very important to have a
sufficiently long follow-up time to capture the events of interest,
and the minimum follow-up time needed depends on both the
aggressiveness of the disease and the type of endpoint} (\cite{Liu2018}, p.401).
See also Othus et al. (\cite{ober:2020}, p.1038),  who ask for
{\it ... Further research ... to identify tests for adequacy of follow-up.}

The test statistics  used in  \cite{Liu2018} are 
\be\label{Q.n}
Q_n =\frac{1}{n} \#\{{\rm uncensored\ observations\ exceeding}\ 2M_u(n)-M(n) \},
\ee
 suggested by Maller and Zhou (\cite{MZ1996}, p.81),  and a similarly constructed  alternative  suggested by Shen \cite{Shen2000}.
Large values of $Q_n$ are associated with sufficient follow-up, thus, provide evidence against  $H_0: \tau_G<\tau_{F}$.
We reject $H_0$ and conclude that follow-up is sufficient if the observed value of $Q_n$ exceeds its 95-th percentile calculated under the null. 

The results in Section \ref{s1} can be used to get expressions for the distribution of $Q_n$.
Conditional on the event
$(M_u(n)=t,M(n)=x)$, the RHS of \eqref{Q.n} is
\be\label{q1}
\frac{1}{n} \sum_{i=1}^n {\bf 1} \{2t-x< T_i, C_i=1\}
=
\frac{1}{n} \sum_{i=1}^n {\bf 1} \{2t-x<T_i^*\le U_i \}
\ee
(with  $Q_n=0$ if $2t-x\ge \tau_H$).
Using Theorem \ref{th1} and Remark (ii) following it,
$nQ_n$ is conditionally distributed as a Binomial $(n, \pi(2t-x))$ rv, where $\pi(2t-x)= \int_{2t-x}^{\tau_H} \Gbar(t)\rmd F^*(t)$ when $0<2t-x< \tau_H$.
Since we know the joint distribution of  $M(n)$ and $M_u(n) $  from Theorem \ref{th4a}, 
we can obtain  formulae for the previously unknown unconditional distribution of $Q_n$, and for similarly constructed statistics. So our present results open up  wide areas of statistical application. 
We leave  further development of   these ideas till later.
%


\section{Proofs for Section \ref{s1}}\label{pfs}
Recall  we assume throughout that $F$ and $G$ are continuous, including at 0 and at their right extremes if these are finite.
Thus the (censored) survival times are all distinct with probability 1. 

\subsection{Proof of Theorem \ref{th1}}\label{s4}
A natural split of the sample into values greater than  or less  than  $ M_u(n) $  is exploited to prove Theorem \ref{th1}.
We need some  notation for sample values less  than  and greater than $ M_u(n) $.
First consider  sample values greater than $ M_u(n) $.
Fix $0<t\le \tau_H$ and let
\be\label{trun}
\big(T_i^>(t), C_i^>(t) \big)_{1\le i\le n}=
\big(T_i^{*,>}(t)\wedge  U_i^>(t),\, {\bf 1}(  T_i^{*,>} (t) \le   U_i^>(t))\big)_{1\le i\le n}
\ee
be a censored sample from  $(  T_i^{*,>}(t),  U_i^>(t))_{i\ge 1}$, where 
 $(  T_i^{*,>}(t))$ and $  (U_i^>(t))$ are two independent sequences of iid  positive random variables whose components  have distributions the same as
\be\label{t2}
T^{*,>}(t) \eqdr \big(T^*|T^*>t\big)
 \ {\rm and}\ 
 U^>(t) \eqdr\big( U|U>t\big).
\ee
(Recall the convention in \eqref{e: conditDistn}).
Analogously, let 
\be\label{cot}
T^{>,c}(t)  \eqdr \big(T^>(t) | U^>(t)<  T^{*,>}(t)\big).
\ee
%



For sample values less  than  $ M_u(n)$,
fix $t>0$ and let
  $(  T_i(t)_{i\ge 1}$ be a sequence of iid  positive random variables 
  having distribution    
\be\label{v2}
T(t) \eqdr \big(T|T<t\big).
\ee

Having set up this preliminary notation we can commence to prove the statement in Theorem \ref{th1}. 
Let $A_k$ and $B_k$ denote arbitrary Borel  sets in $[0,\infty)^k$. 
Take $0<r<n-1$ and $t>0$ and consider
\begin{align}\label{w1}
&
P\big(  M_u(n)>t;\, {\rm there\ are}\ r\, {\rm  censored \ observations\ exceeding}\ 
M_u(n)  \, {\rm and\ they}
  \cr
& \quad {\rm  are\ in}\ M_u(n) + A_r; \,
{\rm  the\ remaining}\ n-r-1\   {\rm  observations\ are\ in}\ B_{n-r-1}\big).
\end{align}
We can calculate this probability as
\bea\label{w2}
&&
 \sum_{\ell=1}^n
\sum_{i_1,\ldots, i_{n-r-1}}   
P\big( C_\ell=1,\, T_\ell> t\vee \max_{ j\in\N_{n-r-1} }   T_{i_j}, \, 
 (T_{i_j})_{ j\in\N_{n-r-1} }  \in  B_{n-r-1}, \cr
&&
(T_i)_{ i\in\NN_n\setminus\{ i_1,\ldots, i_{n-r-1}, \ell \}  }\in T_\ell+A_r\, \& \,  C_i=0,  i\in\NN_n\setminus\{ i_1,\ldots, i_{n-r-1}, \ell \}     \big).  \cr
&&
\eea
In this expression,
$i_1,\ldots, i_{n-r-1}$ are $n-r-1$ unequal integers in $\NN_{n-1}$, distinct from $\ell$,
 and $T_{\ell}$ exceeds all  of the corresponding $T_{i_j}$, and also $T_\ell>t$; 
these $T_{i_j}$ are the observations, both censored and uncensored, smaller than $T_\ell$.
The  remaining $T_i$, of which there are $r$, 
i.e., those with $i\in \NN_n\setminus \{i_1,\ldots, i_{n-r-1}, \ell\}$, are the censored observations exceeding $T_\ell$, which is the largest  uncensored observation in the sample, and those $T_i$ are in $T_\ell+A_r$.

On the event in \eqref{w2}, $T_\ell=T_\ell^*\le U_\ell$
and  $T_{i}=U_{i}<T_i^*$,  for $i\in \NN_n\setminus \{i_1,\ldots, i_{n-r-1}, \ell\}$.
Thus, using   exchangeability to renumber the observations conveniently,
the probability in \eqref{w2} equals 
\bea
&&
n{n-1\choose n-r-1}
P\big(
U_{r+1}\ge  T_{r+1}^*>t\vee \max_{r+2\le i\le n}T_i,\,
(T_i)_{ r+2\le i\le n }  \in  B_{n-r-1}; 
   \cr
&&\cr
&&\hskip1.5cm 
T_i \in  T_{r+1}^*  + A_r,\, T_i^*>U_i> T_{r+1}^*,\,  i\in\N_r
\big).
\eea
Condition on $ T_{r+1}^*$ and integrate to rewrite  this as
\bea\label{w3}
&&n{n-1\choose r} \int_{y=t}^{\tau_H} 
P\big((T_i)_{i\in\N_r} \in  y + A_r, T_i^*>U_i>y, i\in\N_r\big)\cr
&&\cr
&& \qquad
\times 
\Gbar(y)
P\big(y>  \max_ {r+2\le i\le n}T_i, \,
(T_i)_{r+2\le i\le n} \in  B_{n-r-1}\big) \rmd F^*(y),
\eea
where we see that the observations greater than  $ T_{r+1}^*$ are independent 
of those smaller than it, and recall that 
$P(T_{r+1}^*\le y)=F^*(y)$ and $P(U_{r+1}>y)=\Gbar(y)$. 
(Notice that since one or other of $F^*$ or $G$ attributes no mass to values exceeding $\tau_H=\min(\tau_{F^*},\tau_g)$, we can replace $\tau_{F^*}$ and $\tau_G$ by $\tau_H$ in any of the integrals such as appear in \eqref{w3}.)


 The first probability  in \eqref{w3}  is 
\be\label{w11}
P\big((T_i)_{i\in\N_r} \in  y + A_r \big|T_i^*>U_i>y, i\in\N_r\big)
\times  P\big(T_i^*>U_i>y, i\in\N_r\big),
\ee
and  the first probability  in \eqref{w11}  is 
\bea\label{w12}
&&
P\big((T_i^>(y))_{i\in\N_r} \in  y + A_r \big|T_i^{*,>}(y)>U_i^>(y), i\in\N_r\big)\cr
&&=
P\big((T_i^{>,c}(y))_{i\in\N_r} \in  y + A_r\big)
\eea
(using the notation in \eqref{cot}).
The second probability  in \eqref{w11}  is 
\bea\label{w6}
&&
P( {\rm in\ a\ sample\ of\ size}\ r,\, {\rm all\ observations\ are\ censored \ and\ all\ exceed}\ y ) \cr
&&
\hskip1cm 
=\Big(\int_{u=y}^{\tau_H} \Fbar^*(u)\rmd G(u) \Big)^r.
\eea
Consequently  the first probability  in \eqref{w3}   can be written as 
\be\label{w4}
P\big( (T_i^{>,c}(y))_{i\in\N_r} \in y+ A_r\big)
\times  \Big( \int_{z=y}^{\tau_H} \Fbar^* (z)\rmd G(z) \Big)^{r}.
\ee
Using the notation in \eqref{v2}, the second probability in \eqref{w3}    can be written as 
\bea\label{w5}
&&
P\big(   (T_i(y))_{r+2\le i\le n} \in B_{n-r-1}\big)
\times 
 P^{n-r-1} \big(T_i<y \big)\cr
 &&
 =
 P\big(   (T_i(y))_{r+2\le i\le n} \in B_{n-r-1}\big) \times H^{n-r-1}(y).
\eea
So we can write \eqref{w3} and consequently \eqref{w1} for $0<r<n-1$ as
\bea\label{w101}
&&
n{n-1\choose r}
\int_{y=t}^{\tau_H} 
P\big( (T_i^{>,c}(y))_{i\in\N_r} \in y+ A_r\big)
\Big( \int_{z=y}^{\tau_H} \Fbar^* (z)\rmd G(z) \Big)^{r}\cr
&&\cr
&&\hskip1cm  
\times 
P\big(   (T_i(y))_{r+2\le i\le n} \in B_{n-r-1}\big)
 H^{n-r-1}(y)  \Gbar(y) \rmd F^*(y).
\eea

When $r=0$, \eqref{w1} is interpreted as
\bea\label{w111}
P\big(  M_u(n)>t;\,  {\rm  no\ censored \ observation\ exceeds}\  M_u(n), 
  \cr
{\rm  and\ the\ remaining}\ n-1\   {\rm  observations\ are\ in}\ B_{n-1}\big).
\eea
We can calculate this probability as
\bea\label{w222}
&&
 \sum_{\ell=1}^n
P\big( C_\ell=1,\, T_\ell> t\vee 
\max_{ i\in\N_{n}\setminus\{\ell\} } T_i,\,
 (T_{i})_{ i\in\N_{n}\setminus\{\ell\}  }  \in  B_{n-1}\big) \cr
&&=
n  \int_{y=t}^{\tau_H} 
P\big(   (T_i(y))_{2\le i\le n} \in B_{n-1}\big)
 H^{n-1}(y)  \Gbar(y) \rmd F^*(y),
\eea
which agrees with \eqref{w101}  when $r=0$ if we interpret $y+A_0=\emptyset$, the empty set.

When $r=n-1$, \eqref{w1} is interpreted as
\bean
&&
P\big(  M_u(n)>t,\, {\rm  the\ remaining}\ n-1\   {\rm  observations\ are\ censored} \cr
&& \qquad
{\rm  and \ in}\ M_u(n)+A_{n-1}\big)\cr
&&=
 \sum_{\ell=1}^n
P\big( C_\ell=1,\, T_\ell> t,
U_i< T_i^*,  i\in\N_{n}\setminus\{\ell\},
 (U_{i})_{ i\in\N_{n}\setminus\{\ell\}  }  \in  T_\ell+A_{n-1}\big) \cr
&&=
n  \int_{y=t}^{\tau_H} 
P\big(   (U_i(y))_{2\le i\le n} \in y+A_{n-1}\big)
\Big( \int_{z=y}^{\tau_H} \Fbar^* (z)\rmd G(z) \Big)^{n-1}
  \Gbar(y) \rmd F^*(y),
\eean
which agrees with \eqref{w101}  when $r=n-1$ if we interpret $B_{0}=\emptyset$.

Thus we see that \eqref{w101} holds for all $0\le r\le n-1$
with the appropriate interpretations.
Choosing $A_r=[0,\infty)^r$ and $B_{n-r-1}=[0,\infty)^{n-r-1}$ in  \eqref{w101} 
and taking complements gives \eqref{j0}, 
and from \eqref{j0} we have
\bea\label{w10a}
&&
P\big(M_u(n)\in \rmd t, N_c^>(M_u(n))=r\big)\cr
&&\cr
&& =n{n-1\choose r}
\Big( \int_{z=t}^{\tau_G} \Fbar^* (z)\rmd G(z) \Big)^{r}
H^{n-r-1}(t)   \Gbar(t)\rmd F^*(t),
\eea
for $0\le r\le n-1$. 
So \eqref{w101} and consequently \eqref{w1}  can be rewritten as
\bea\label{w102}
&&
 \int_{y=t}^{\tau_H} 
P\big( (T_i^{>,c}(y))_{i\in\N_r} \in y+A_r\big)
\times 
P\big(   (T_i(y))_{r+2\le i\le n} \in B_{n-r-1}\big)\cr 
&&\hskip4cm  
\times
P\big(M_u(n)\in \rmd y, N_c^>(M_u(n))=r\big).
\eea
Recalling the definitions in \eqref{trun}--\eqref{v2}, this gives \eqref{w10b}.
 \halmos

\medskip\noindent{\bf Remarks.}\
(i)\ 
Conditional on $\{M_u(n)=t, N_c^>(n)=r\}$, $t>0$, $0<r<n-1$, 
 the independent components into which the sample splits are
\be\label{is}
\{(T_i,C_i): T_i<t, 1\le i\le n-r-1\}
\ {\rm and}\
\{(T_i,C_i): t<U_i< T_i^*, 1\le i\le r\},
\ee
where $T_i=T_i^* \wedge U_i$ and $C_i={\bf 1}_{\{T_i^*\le U_i\}}$, $1\le i\le n$.
The cases $r=0$ and $r=n-1$ are  included in  Theorem \ref{th1} with appropriate interpretations, as previously noted.

(ii)\ 
We can  include the case $r=n$ in  Theorem \ref{th1} by deriving  the conditional distribution of the observations given $M_u(n)=0$.
When $r=n$ 
 all observations are censored 
 in which case the second component  in \eqref{is} contains the whole sample and the first component  in \eqref{is} is empty.  By convention we then set $ M_u(n)=0$.
Then  we can calculate
\bea\label{c2a}
P\big(T_i\le t_i, 1\le i\le n,M_u(n)=0\big)
&=&
 P\big(T_i^*>U_i, U_i\le t_i, 1\le i\le n\big)\cr
&=&
\prod_{i=1}^n\int_{z=0}^{t_i}  \Fbar^* (z)\rmd G(z),
\eea
and thus, for $0<x\le \tau_H$, 
\be\label{c1a}
P\big(M_u(n)=0, M(n)\le x\big)
= 
\Big(\int_{z=0}^{x}  \Fbar^* (z)\rmd G(z)\Big)^n.
\ee
Dividing \eqref{c2a} by \eqref{c1a} with $x=\tau_H$ gives the required conditional distribution as:
\be\label{c0a}
P\big(T_i\le t_i, 1\le i\le n|M_u(n)=0\big)
= 
\frac{\prod_{i=1}^n\int_{z=0}^{t_i}  \Fbar^* (z)\rmd G(z)}
{\big(\int_{z=0}^{\tau_H}  \Fbar^* (z)\rmd G(z)\big)^n},\ t_i\ge 0.
\ee

The next corollary is required for Remark (ii) following  Theorem \ref{th1}.
Define integers
$I_n^<:= \{i\in\N_n: T_i<M_u(n)\}$
and 
$I_n^>:= \{i\in\N_n: T_i>M_u(n)\}$,
and let
$\sigma^<$ be the smallest  $\sigma$-field 
making $(T_i,C_i)_{i\in I_n^<}$ measurable;
and likewise let 
$\sigma^>$ be the smallest  $\sigma$-field 
making $(T_i,C_i)_{i\in I_n^>}$ measurable.

\begin{corollary}\label{c1}
Let 
$A^<$ be any event in $\sigma^<$ and
$A^>$ any event in $\sigma^>$.
Then for any Borel $B\subseteq [0,\infty)$, $t>0$, $0\le r\le n-1$, 
\bea\label{cor0}
&&
P\big(A^<\big| A^>, M_u(n)=t, M(n)\in B, N_c^>(M_u(n))=r)\cr
&&\cr
&& \qquad
=P\big(A^<\big|M_u(n)=t, N_c^>(M_u(n))=r).
\eea
\end{corollary}

\noindent{\bf Proof of Corollary \ref{c1}}
For a Borel $B_1\subseteq\R_+= [0,\infty)$, $t>0$, $0\le r\le n-1$,
we have 
\begin{align}
P\big(A^<, &A^>, M_u(n)\in B_1, M(n)\in B, N_c^>(M_u(n))=r\big) \label{e:cor1}\\
=& \int_{t\in B_1} 
P\big(A^<, A^>, t\vee \max_{i\in I_n^>} U_i\in B
   \big|M_u(n)=t, N_c^>(M_u(n))=r\big)\nonumber\\
&\qquad \times P(M_u(n)\in\rmd t, N_c^>(M_u(n))=r) .\nonumber\\
\intertext{When $N_c^>(M_u(n))=r$, an event $A^<$ in $
  \sigma^<$ is of the form  
$A^<= \big\{I_n^<= \{ i_1, \ldots,i_{n-r-1}\}, \,
 (T_i,C_i)_{i\in I_n^<} \in  B_{n-r-1} \big\}$, 
where 
$i_1, \ldots,i_{n-r-1}$ are unequal integers in $\N_n$ 
and $ B_{n-r-1} $ is Borel in $(\R_+ \times \{0,1\})^{n-r-1}$.
Then, on the event $\{M_u(n)=t,\, N^>(M_u(n))=r\}$, 
 we obtain $A^<(t)$ by writing $T_i(t)$ for $T_i$ in $A^<$.
Similarly, we can formulate  $A^>$ and $A^>(t)$. Using these remarks,
  \eqref{e:cor1} becomes}
=&\int_{t\in B_1} P(A^<(t))\,\times \, P\big( A^>(t), t\vee \max_{i\in I_n^>} U_i^>\in B \big) \nonumber\\
 &\qquad \times P(M_u(n)\in\rmd t, N_c^>(M_u(n))=r), \label{e:rhs}\end{align}
where the factorisation of the integrand is justified by  Theorem \ref{th1}.
Now for $t>0$
\bea\label{cor2}
&&
\int_{0<y\le t} 
P\big(A^>, M_u(n)\in\rmd y, M(n)\in B, N_c^>(M_u(n))=r)\cr
&&
=
\int_{0<y\le t}  
P\big( A^>(y), y\vee \max_{i\in I_n^>} U_i^>\in B \big)\cr
&&\hskip4cm 
\times P(M_u(n)\in \rmd y, N_c^>(M_u(n))=r),
\eea
by   Theorem \ref{th1} again.
From this we see that 
\eqref{e:rhs}
equals
\be\label{cor3}
\int_{t\in B_1} 
P(A^<(t))\,
\times \,
P\big( A^>, M_u(n)\in \rmd t, M(n) \in B, N_c^>(M_u(n))=r).
\ee
Comparing the LHS of \eqref{e:rhs} with \eqref{cor3} shows that
\bea\label{cor4}
&&
P\big(A^<\big| A^>, M_u(n)=t, M(n)\in B, N_c^>(M_u(n))=r\big)\cr
&& 
=
P(A^<(t))
= 
P\big(A^<\big| M_u(n)=t, N_c^>(M_u(n))=r\big),
\eea
which is \eqref{cor0}.  \halmos

\subsection{\bf Proof of Theorem \ref{th4a}}\label{pfth4a}
\noindent{\bf  Part  (i)}\ 
 Keep   $0\le x\le \tau_H$,  $0\le t\le \tau_H$ and $0\le r\le n-1$, and calculate, using  Theorem \ref{th1},
\begin{align}\label{j7}
&
P\big(0<M(n) \le x\big| M_u(n)=t,\, N_c^>(M_u(n))=r\big)\cr
&
=
 {\bf 1}_{\{t\le x\}} 
P\big(
\max_{1\le i\le n}T_i \le x\big| M_u(n)=t,\, N_c^>(M_u(n))=r\big) \cr
&
=
 {\bf 1}_{\{t\le x\}} 
P\big(T_i(t)\le x, 1\le i\le n-r-1\big) 
P\big(T_i^{>,c}(t)\le x, 1\le i\le r\big), \ \
\end{align} 
where if $r=0$ the second probability on the RHS is taken as 1 and if $r=n-1$ the first probability on the RHS is taken as 1.
The first probability on the RHS of \eqref{j7} also equals 1 when $t\le x$. So,
recalling the definition of $T_1^{>,c}(t)$ in \eqref{cot}, we get
\begin{align}\label{j8}
&P\big(0<M(n) \le x\big| M_u(n)=t,\, N_c^>(M_u(n))=r\big)\cr
&=
{\bf 1}_{\{t\le x\}} 
P^r\big(T_1^{>,c}(t)\le x\big)
=
{\bf 1}_{\{t\le x\}} 
\left(
\frac{\int_{z=t}^{x}  \Fbar^* (z)\rmd G(z)}
{\int_{z=t}^{\tau_H}  \Fbar^* (z)\rmd G(z)}\right)^r.
\end{align} 
 Next keep   $0<t \le x\le \tau_H$ and $0\le r\le n-1$, and use
 this together with  \eqref{j0} 
 to calculate
\bea\label{j9}
&&P\big(0< M_u(n)\le t, 0\le  M(n)\le x, N_c^>(M_u(n))=r\big)\cr
&&=
\int_{y=0}^t P\big(0\le M(n)\le x|M_u(n)=y, N_c^>(M_u(n))=r\big)\cr
&&
\hskip6cm 
\times 
P\big(M_u(n)\in \rmd y, N_c^>(M_u(n))=r\big)\cr
&&=
\int_{y=0}^t 
{\bf 1}_{\{y\le x\}} \frac{\Big(\int_{z=y}^{x}  \Fbar^* (z)\rmd G(z)\Big)^r}
{\Big(\int_{z=y}^{\tau_H}  \Fbar^* (z)\rmd G(z)\Big)^r}\cr
&&\hskip2.5cm 
\times  n{n-1\choose r}
\Big( \int_{z=y}^{\tau_H} \Fbar^* (z)\rmd G(z) \Big)^{r}H^{n-r-1}(y) \Gbar (y)\rmd F^*(y)\cr
&&\hskip0.5cm  = n{n-1\choose r}
\int_{y=0}^t 
\Big(\int_{z=y}^{x}  \Fbar^* (z)\rmd G(z)\Big)^r   H^{n-r-1}(y) \Gbar (y)\rmd F^*(y).
\eea
(Observe that ${n\choose r}(n-r) = n{n-1\choose r}$.)

Add over $0\le r\le n-1$ in \eqref{j9} and recall \eqref{F^*GH}  to get
\bea\label{j6}
&&
P\big(0<M_u(n)\le t,\, 0\le M(n)  \le x\big) \cr
&&
=n \int_{y=0}^ {t}
\Big( \int_{z=y}^{x} \Fbar^* (z)\rmd G(z) +H(y)\Big)^{n-1}  \Gbar(y)\rmd F^*(y)\cr
&&\cr
&&=
\Big( \int_{z=t}^{x} \Fbar^* (z)\rmd G(z) +H(t)\Big)^{n} 
 - \big(\int_{z=0}^{x} \Fbar^*(z)\rmd G(z)\big)^n.
\eea
Adding in the value for $t=0$ in \eqref{c1a}  gives
\ben 
P\big(0\le M_u(n)\le t,\, 0<M(n)  \le x\big) 
=
\Big( \int_{z=t}^{x} \Fbar^* (z)\rmd G(z) +H(t)\Big)^{n}
\een
 for $0\le t \le x\le \tau_H$, 
and hence  the third line  on the RHS of  \eqref{JD}. 

For $0\le x<t \le \tau_H$, take $t=x$ in both sides of \eqref{j6} to get
\ben 
P\big(0< M_u(n)\le x, 0\le  M(n)\le x\big)
=
H^n(x) 
 - \big(\int_{z=0}^{x} \Fbar^*(z)\rmd G(z)\big)^n.
\een
Adding in the value for $t=0$ in \eqref{c1a}  gives
\ben 
P\big( 0\le  M(n)\le x\big)=
P\big(0\le M_u(n)\le x,\,  0\le  M(n)\le x\big)
=
H^n(x)
\een
(as it should), 
and hence  the second line  on the RHS of  \eqref{JD}. 

Setting all $t_j=t\ge 0$ in  \eqref{c2a} 
gives
\be\label{j14}
P\big(M_u(n)=0, M(n)\le x\big)  
= \Big(\int_{z=0}^{x} \Fbar^* (z)\rmd G(z)\Big)^{n},
\ee
which is  the first line in \eqref{JD}. 

\noindent {\bf Part  (ii)}\ 
Taking $x=\tau_H$ in  the first and third  lines  on the RHS of   \eqref{JD} gives
\ben 
P\big(0\le M_u(n) \le  t\big)= 
\Big( \int_{z=t}^{\tau_H} \Fbar^* (z)\rmd G(z)+H(t) \Big)^{n},
\een
which is \eqref{4f} in terms of the righthand formula in \eqref{4g}.
 The lefthand formula in \eqref{4g} comes from an integration by parts.  
For the right extreme $\tau_J$ of the distribution $J$ we have $\tau_J=\tau_{F}\wedge \tau_G$.
This is established by checking the behaviour of the second integral in \eqref{4g} in the cases: Case 1, $p=1$, and (a) $\tau_{F^*}=\tau_{F}\le \tau_G$, or else (b) $\tau_{F^*}=\tau_{F}>\tau_G$;
and Case 2, $0<p<1$, in which case $\tau_{F^*}=\infty$, and again we may have, (a) $\tau_{F}\le \tau_G$, or else (b) $\tau_{F}>\tau_G$.
 \halmos

%
%
%

\subsection{Proof of Theorem \ref{th4}}\label{pth4}
To prove  \eqref{12a}, take $0<u\le \tau_H$ and write
\bea\label{14a}
&&P\big( M(n)  - M_u(n)  \le u \big)\cr
&&=
P\big(M(n) =M_u(n)  \big) +P\big(0< M(n)  - M_u(n)
  \le u \big) 
=:
   A+B. \
\eea
Decompose the component $A$ as
\bean 
&&A=
P\big( M(n) =M_u(n)  \big) \cr
&&\cr
&&
=
P\big(M_u(n)=0,  M(n) =M_u(n)  \big) 
+P\big(0< M(n) =M_u(n)  \big)=:
 A_1+A_2.
\eean
Here $A_1= P\big(M(n)=0 \big)=0$. 
For $A_2$ we calculate 
\bea\label{j5}
&&
P\big(0< M(n) =M_u(n) \le t \big) 
=
P\big(0< M_u(n)\le t,  N_c^>(M_u(n))=0\big)\cr
&&=
\sum_{\ell=1}^n  
P\big( 0<T_\ell^*\le t\wedge U_\ell, 
 T_i\le  T_\ell^*, i\in \NN_n, i\ne \ell \big)\cr
 &&\cr
 &&=
n \int_{y=0}^{ t}    H^{n-1}(y) \Gbar(y)\rmd F^*(y).
\eea
Taking $t=\tau_H$ in \eqref{j5} we find
\be\label{star}
A_2
=P\big(0< M_u(n)= M(n)   \le \tau_H \big) 
= n\int_{y=0}^{\tau_H} H^{n-1}(y)\Gbar(y) \rmd F^*(y).
\ee
Further decompose $B$ in \eqref{14a} as
\bean 
&&P\big(0< M(n)  - M_u(n)  \le u \big) 
=
P\big( M_u(n)=0, 0< M(n)   \le u \big) \cr
&&
+
 P\big( M_u(n)>0, 0< M(n)  - M_u(n)  \le u \big) 
 =:
 B_1+B_2.
\eean
By \eqref{j14}, $B_1= \big(\int_{z=0}^{u} \Fbar^*(z) \rmd G(z)\big)^n$.
To get $B_2$, we calculate, for $0<u\le \tau_H$, 
\bean 
&&
B_2= 
 P\big( M_u(n)>0, 0< M(n)  - M_u(n)  \le u \big) \cr
&&=
\int_{0<t\le \tau_H} P\big(0\le M(n) \le t+u\big| M_u(n)=t \big)P\big(M_u(n)\in \rmd t\big)\cr
&&=
\int_{0<t\le \tau_H}
\frac{\Big(\int_{z=t}^{\min((t+u), \tau_H)} \Fbar^*(z) \rmd G(z) +H(t)\Big)^ {n-1} }
{\Big(\int_{z=t}^{\tau_H} \Fbar^*(z) \rmd G(z) +H(t)\Big)^ {n-1}}\cr
&& \hskip 2cm
\times
n \Big(\int_{z=t}^{\tau_H} \Fbar^*(z) \rmd G(z) +H(t)\Big)^ {n-1} \Gbar(t) \rmd F^*(t) \cr
&&\cr
&&
=  n\int_{0<t\le \tau_H}
\Big(\int_{z=t}^{\min((t+u), \tau_H)} \Fbar^*(z) \rmd G(z) +H(t)\Big)^ {n-1} 
 \Gbar(t) \rmd F^*(t). \cr
&&
\eean
The conditional distribution of $M(n)$ given $M_u(n)$ follows from \eqref{4f} and \eqref{j6}.
For $P\big( M(n)  - M_u(n)  \le u \big)$ we add $A+B= A_1+A_2+B_1+B_2 = A_2+B_1+B_2$.
But note that $A_2$  is the same as $B_2$ with $u$ set equal to 0, so that $P\big( M(n)  - M_u(n)  \le u \big)$ is given by $B_1+B_2$ 
(with $u= 0$ allowed  in $B_2$). This verifies \eqref{12a}. 

To prove  \eqref{12b}, take $v\ge 1$ and write
\begin{align*} 
&P\big( M_u(n)  >0,  M(n)  \le v M_u(n)  \big) \cr
&=
P\big(0<M(n) =M_u(n)  \big) +P\big( 0< M_u(n) < M(n)  \le v M_u(n)  \big) \cr
&=: 
A+C.
\end{align*}
Here 
$A=A_1+A_2=A_2$ as in \eqref{star}.
When $v=1$, $C=0$. 
When $v>1$, $C$ is
\bea\label{16v}
&&
 P\big(0< M_u(n)< M(n)  \le v  M_u(n)  \big) \cr
&&=
\int_{0<t\le \tau_H} P\big(t< M(n) \le vt\big| M_u(n)=t \big)P\big(M_u(n)\in \rmd t\big)\cr
&&=
\int_{0<t\le \tau_H}
\frac{\Big(\int_{z=t}^{\min(tv, \tau_H)} \Fbar^*(z) \rmd G(z) +H(t)\Big)^ {n-1} }
{\Big(\int_{z=t}^{\tau_H} \Fbar^*(z) \rmd G(z) +H(t)\Big)^ {n-1}}\cr
&& \hskip 4cm
\times
n \Big(\int_{z=t}^{\tau_H} \Fbar^*(z) \rmd G(z) +H(t)\Big)^ {n-1} \Gbar(t) \rmd F^*(t) \cr
&&\cr
&&
=  n\int_{0<t\le \tau_H}
\Big(\int_{z=t}^{\min(tv, \tau_H)} \Fbar^*(z) \rmd G(z) +H(t)\Big)^ {n-1} 
 \Gbar(t) \rmd F^*(t). 
\eea
For $P\big( M_u(n)>0, M(n) \le v M_u(n)  \big)$ we take $A+C= A_2+C$.
But $A_2$  is the same as $C$ with $v$ set equal to 1.
So $P\big( M_u(n)>0,  M(n) \le vM_u(n)  \big)$ is given by $C$ 
(with $v=1$ allowed  in $C$). Then to obtain \eqref{12b} divide \eqref{16v} as it is by \eqref{16v} with $v$ set equal to $\infty$.
 \halmos

\subsection{Proof of Theorem \ref{th2}}\label{th2.5}

Here we use the d\'ecoupage de  L\'evy (e.g., Resnick \cite{resnickbook:2008}, p.212), which we used in \cite{MR2020} and which we now review.
Relative to a sequence $\{(T_i^*,U_i), i \geq 1\}$, as specified in Section \ref{s1}, define random
indices $K_i^u$ and  $K_i^c$ by
\begin{align}\label{e:Kgreater}
  \Kless_0=&0, \quad \Kless_i =\inf \{m>\Kless_{i-1} :
   T_m^*\leq U_m\}, \ i\ge 1, \cr 
   \ {\rm and} \quad
\Kgreater_0=&0, \quad \Kgreater_i =\inf \{m>\Kgreater_{i-1} : T_m^*> U_m\}, \ i\ge 1.
\end{align}
The sequences $\{T_{\Kless_i}, i\ge 1\}$ and 
 $\{T_{\Kgreater_i}, i\ge 1\}$, select out 
the subsequences of uncensored and censored  observations, respectively. 
Note that  $T_{\Kless_i}^*=  T_{\Kless_i}$ and
$U_{K_i^c}=T_{K_i^c}$. 
The $T_{\Kless_i}$ and $T_{\Kgreater_i}$ are iid with respective distributions
\be\label{tdis}
P(T_{\Kless_i}\in A) =P(T_1\in A|T_1^*\le U_1)
\    {\rm and} \
P(T_{\Kgreater_i}\in A) =P(U_1\in A|T_1^*>U_1),
\ee
 for Borel $A \subseteq[0,\infty)$. 
%
Both  subsequences $\{(T_{\Kless_i}^*,U_{\Kless_i} )\}$ and $\{(T_{\Kgreater_i}^*,U_{\Kgreater_i}) \}$ are comprised of iid  random vectors.
Furthermore, the three sequences
\be\label{e:indep}
\{(T_{\Kless_i}^*,\, U_{\Kless_i} ),\, i\geq
1\},\, \{(T_{\Kgreater_i}^*,\, U_{\Kgreater_i}), \, i\geq 1\}, \,
\{N_u(i): =\sum_{m=1}^i {\bf 1}_{\{ T_m^*\le U_m\}},\, i\geq 1\},
\ee 
 are independent of each other. 
 
 Recall the notations for the numbers $N_u(n)$, $N_c(n)$,  $N_u^<(M_u(n))$, $N_c^<(M_u(n))$ and $N_c^>(M_u(n))$  in \eqref{Nd1}--\eqref{Nd6} for a sample of size $n$.
 In calculating \eqref{w15a}, we have a sample of size $n$ and  condition on $M_u(n)=t$ with $t>0$.
 This means that there is at least once uncensored observation,
 so $N_u(n)\ge 1$.     
We index the $N_u(n)$ uncensored observations as
 $(T_{\Kless_i}^*)_{1\le i\le N_u(n)}$ 
 and the   $N_c(n)= n-N_u(n)-1$ censored observations as 
  $(U_{\Kgreater_i})_{1\le i\le N_c(n)}$. 


\noindent
{\bf Part (i)} \
To prove \eqref{w15a} we begin  by calculating, for nonnegative integers $r, s,k$, 
with $0\le r,s,k\le n-1$, $r+s+k=n-1$, and $t>0$, 
the probability
\be\label{cd1}
P\big( N_c^>(M_u(n)) =r,\,  N_c^<(M_u(n)) =s
\big|  M_u(n)=t,\,  N_u(n)=k+1\big).
\ee
 Substituting for the definitions of the numbers in \eqref{Nd1}--\eqref{Nd6} 
 we have to calculate  for \eqref{cd1} the probability of the event
\bea\label{cd2}
\Big\{
\sum_{i=1}^{n-k-1} {\bf 1}_{\{ t<  U_{\Kgreater_i} < T_{\Kgreater_i}^*\}}=r \Big\}, \ 0\le r\le n-1
 \eea
(note that the requirement $N_c^<(M_u(n)) =s$ in \eqref{cd1} is redundant because we must have $  N_c^<(M_u(n))=n-r-k-1$),
 conditional on the event
 \be\label{cd22}
 \Big\{ \max_{1\le i\le k+1} T_{\Kless_i}^*=t,\,  
\sum_{i=1}^n {\bf 1}_{\{ U_{\Kless_i} \ge  T_{\Kless_i}^*\}} =k+1\Big\}, \ 0\le k\le n-1.
\ee

Denote  the sum in \eqref{cd2} by 
\be\label{cd33}
N_c^>(t,n-k-1):= \sum_{i=1}^{n-k-1} 
{\bf 1}_{\{ t<  U_{\Kgreater_i} < T_{\Kgreater_i}^*\}}.
\ee
%
As a result of the d\'ecoupage in \eqref{e:indep}, the variables in
 \eqref{cd22}, having index ``$u$", 
 are independent of those in \eqref{cd2}, 
having index ``$c$". So the conditioning event in \eqref{cd22} is independent of the  event in \eqref{cd2}.
Thus the probability in \eqref{cd1} equals
$P(N_c^>(t,n-k-1)=r)$. 
With this notation we can write what we have  proved so far as
\bea\label{cd4}
&&
P\big( N_c^>(M_u(n)) =r,\,  N_c^<(M_u(n)) =s
\big|  M_u(n)=t,\,  N_u(n)=k+1\big)\cr
 &&
  \hskip4cm 
   = P\big(N_c^>(t,n-k-1)=r\big).
 \eea
%
Rewriting the conditional probability in \eqref{cd1} using \eqref{cd4}, and the fact that $N_u^<(M_u(n))=N_u(n)-1$ on $\{M_u(n)>0\}$, 
 we obtain
\bea\label{cd5}
&&
P\big( N_c^>(M_u(n)) =r,\,  N_c^<(M_u(n)) =s,\,  N_u^<(M_u(n))=k|  M_u(n)=t\big)  \cr
&&\cr
 &&
=P\big(N_c^>(M_u(n))=r,\, N_c^<(M_u(n))=s| M_u(n)=t, \, N_u(n)=k+1 \big)\cr
&& \hskip6cm 
\times 
 P\big(N_u(n)=k+1|  M_u(n)=t\big)\cr
 &&\cr
 &&
 = P\big(N_c^>(t,n-k-1)=r\big)
 \times 
 P\big(N_u^<(M_u(n))=k|  M_u(n)=t\big).
 \eea

 The first probability on the RHS of \eqref{cd5} is by \eqref{cd33} a binomial with success probability
 \bean
P(t<  U_{\Kgreater_i} < T_{\Kgreater_i}^*)
&=&
\frac{P(t<U_1<T_1^*)} { P(U_1<T_1^*)}
= \frac{\int_{y=t}^{\tau_H} \Fbar^* (y)\rmd G(y)}
{ \int_{y=0}^{\tau_H} \Fbar^* (y)\rmd G(y)}\cr 
&=&
 \frac{\int_{y=t}^{\tau_H} \Fbar^* (y)\rmd G(y)}
{ \int_{y=t}^{\tau_H} \Fbar^* (y)\rmd G(y)+H(t) }
\times \frac {J(t)}{p_c},
\eean
 where $J(t)$ is  defined in \eqref{4g}
 and $ p_c:=  \int_{y=0}^{\tau_H} \Fbar^* (y)\rmd G(y)$.
 Recalling the notation in \eqref{pdefs}, we can thus write 
 \ben
 P(t<  U_{\Kgreater_i} < T_{\Kgreater_i}^*)
 =
 p_c^>(t)  \times \frac {J(t)}{p_c}.
 \een
  A similar computation gives
  \ben
1- p_c^>(t)  \times \frac {J(t)}{p_c}
 =
 p_c^<(t)  \times \frac {J(t)}{p_c}.
 \een
 Thus the first probability on the RHS of \eqref{cd5} equals
 \be\label{548a}
 {n-k -1\choose r} (p_c^>(t))^r
(p_c^<(t))^{n-k-1-r}  \times \Big(\frac {J(t)}{p_c}\Big)^{n-k-1}.
 \ee
 
  For the second probability on the RHS of \eqref{cd5} we have to carry out a computation similar to that in the proof of Theorem \ref{th1}. We find after some calculation 
  \bea\label{cd14}
  &&
P\big(N_u^<(M_u(n))=k, \,0< M_u(n)\le t\big)\cr
&& \hskip2cm 
= {n \choose k+1}
 \Big(\int_{y=0}^{t} \Gbar (y)\rmd F^*(y)\Big)^{k+1}
  \times p_c^{n-k-1}.
\eea
This is valid for $0\le k\le n-1$.

As a check on \eqref{cd14}, adding \eqref{cd14}  over $0\le k\le n-1$ gives the expression in \eqref{4f} for $P\big(M_u(n)\le t\big)$,
and differentiating  \eqref{4f} gives a formula for 
 $P\big(M_u(n)\in \rmd t\big)$. 
Now 
divide the formula for 
 $P\big(M_u(n)\in \rmd t\big)$ into the corresponding  differential of \eqref{cd14}. (Formally, we calculate 
 Radon-Nikodym derivatives.)
 Then, recalling \eqref{pdefs}, we arrive at
  \be\label{cd15}
P\big(N_u^<(M_u(n))=k|M_u(n)= t\big)
= {n-1 \choose k} (p_u^<(t))^{k} 
 \times \Big(\frac {p_c}{J(t)}\Big)^{n-k-1}.
\ee
Multiplying \eqref{548a} and \eqref{cd15} together and setting $s=n-r-k-1$ gives \eqref{w15a}.

\noindent
{\bf Part (ii)}\
The binomial distributions are immediate from \eqref{w15a}.

\noindent
{\bf Part (iii)}\ This also follows directly from \eqref{w15a} 
by a convolution calculation. 

\noindent
{\bf Part (iv)}\ 
For this we use  the identity
\bea\label{w24}
&&
P\big(N_c^>(M_u(n))=r \big|N_c(n)=\ell, M_u(n)=t  \big) \cr
&&\cr
&&
=
\frac{P\big(N_c^>(M_u(n))=r, N_c(n)=\ell\big| M_u(n)=t  \big)}
{P\big(N_c(n)=\ell  \big| M_u(n)=t  \big)}.
\eea
Here the numerator equals 
\ben
P\big(N_c^>(M_u(n))=r, N_c^<(n)=\ell-r,
N_u(M(n))=n-1-\ell\big| M_u(n)=t  \big)
\een
for which we can obtain a formula from \eqref{w15a}, 
and the denominator is given by the binomial distribution in Part (iii). 
Then it's easily verified that \eqref{w24} is a binomial probability 
with $p_c^+(t)$ defined as in \eqref{pc+}.    \halmos

%

%



\section{Proof of Theorem \ref{FGu}}\label{sspa}
Take any sequences $a_n>0$, $b_n>0$,  $a_n\to\infty$, $b_n\to\infty$, and use the identity $P(A\cap B)= 1+P(A^c\cap B^c)-P(A^c)-P(B^c)$ to write, for any $u,v>0$,  
\bea\label{ap1}
&&
P\big( a_n(\tau_H- M(n))\le u,\ b_n(\tau_J- M_u(n))\le v \big) \cr
&&= 
P\big(M(n)\ge \tau_H-u/a_n,\  M_u(n)\ge \tau_J-v/b_n \big) \cr
&&= 
1+
P\big(M(n)< \tau_H-u/a_n,\  M_u(n)<\tau_J-v/b_n \big)\cr
&&\hskip1.5cm 
- P\big(M(n)< \tau_H-u/a_n \big)
-P\big( M_u(n)< \tau_J-v/b_n \big).
\eea
Recall we assume throughout that $F$ and $G$ are continuous distributions.

\smallskip
\noindent {\bf Case 1:}\  Assume $\tau_{F}<\tau_G<\infty$ and $p<1$.
Then  $\tau_J=\tau_{F}<\tau_H=\tau_G<\tau_{F^*}=\infty$.
Assume also \eqref{F^*GD}. 
For $n\in \NN$ define
$a_n:= \sup\{x>1/\tau_G: \Gbar(\tau_G-1/x)\ge 1/n\}$.
Then $a_n\uto\infty$, $n\Gbar(\tau_G-1/a_n)\to 1$, and by the first relation in \eqref{F^*GD}, 
 $a_n\sim (na_GL_G(1/a_n)) ^{1/\gamma} $ as $n\to\infty$.
Similarly we can choose $b_n$ to satisfy  $n\Fbar(\tau_F-1/b_n)\to 1$, and then  $b_n\sim (na_GL_F(1/b_n)) ^{1/\beta}\to \infty $ as $n\to\infty$.
 
Since $\tau_H>\tau_J$ we can assume $a_n$ and $b_n$ are large enough for 
$\tau_H-u/a_n>\tau_J=\tau_{F} >\tau_J-v/b_n= \tau_{F}-v/b_n$
for any $u,v>0$. 
Recognizing this we can use \eqref{JD}, \eqref{4f}, and \eqref{4h} to express the RHS of \eqref{ap1} as
\bea\label{ap2}
&&1+ \Big(\int_{\tau_J-v/b_n}^{\tau_H-u/a_n}  \Fbar^* (z)\rmd G(z)+ H( \tau_J-v/b_n) \Big)^n  \cr
&&\hskip2cm 
-H^n( \tau_H-u/a_n)  -J^n( \tau_J-v/b_n).
\eea
Recall $\Fbar^*(z)=1-pF(z)$ and $\Hbar(z)=\Fbar^*(z)\Gbar(z)$.
Hence by  \eqref{F^*GD}, the relation  $a_n^\gamma\sim na_GL_G(1/a_n) $,
and  the slow variation of $L_G$,
\bea\label{ap3}  
n\Hbar( \tau_H-u/a_n) 
&=&
n\big(1-p F( \tau_G-u/a_n)\big)\Gbar( \tau_G-u/a_n)\cr
&\sim &
 (1-p) a_G \frac{ u^\gamma  L_G(u/a_n) }{  L_G(1/a_n) }
 \to (1-p)u^\gamma, \ n\to\infty.
\eea
Here $ F( \tau_G-u/a_n)=1$ since  $\tau_G-u/a_n>\tau_{F}$. 
Then by a  standard approximation 
\be\label{ap4}  
\lim_{n\to\infty} H^n( \tau_H-u/a_n) 
= e^{-\lim_{n\to\infty} n\Hbar( \tau_H-u/a_n)} = e^{-(1-p)u^\gamma }, \ u>0.
\ee

For the $J$ term in \eqref{ap2} use \eqref{4g}, $\tau_J=\tau_{F}<\tau_H=\tau_G$  and the mean value theorem for integrals to write
\be\label{ap5}  
n\Jbar( \tau_J-v/b_n) 
=
np \int_{\tau_{F}-v/b_n}^{\tau_{F}}  \Gbar (z)\rmd F(z)
=
np\Gbar(\tau_{F}-z_n)  \Fbar( \tau_{F}-v/b_n)
\ee
where $0 \le z_n \le v/b_n$.
The slow variation of $L_{F}$ and
 $b_n\sim (na_FL_F(1/b_n)) ^{1/\beta}$ imply
$\lim_{n\to\infty} n \Fbar( \tau_{F}-v/b_n)  = v^\beta$.
Since $\Gbar(\tau_{F}-z_n)\to \Gbar(\tau_{F})$,
the RHS of \eqref{ap5}  has limit $p\Gbar(\tau_{F})v^\beta $.
Thus
\be\label{ap6}  
\lim_{n\to\infty} 
J^n( \tau_J-v/b_n)  = e^{- p\Gbar(\tau_{F})v^\beta }.
\ee

Now for the integral term in \eqref{ap2}, subtract from 1 the expression in parentheses,  
recall that $\tau_H=\tau_G$,  $\tau_J=\tau_{F}$ 
and 
 $\tau_H-u/a_n>\tau_J=\tau_{F}> \tau_J-v/b_n$, and calculate as follows:
\bea\label{ap5a}
&&
1-  \int_{\tau_J-v/b_n}^{\tau_H-u/a_n}\Fbar^*(z) \rmd G(z) 
-H( \tau_J-v/b_n)\cr
&&\cr
&&=
\Hbar( \tau_{F}-v/b_n)-  \int_{\tau_{F}-v/b_n}^{\tau_G-u/a_n}   \big(1-p+p \Fbar (z)\big)\rmd G(z)\cr
&&\cr
&&=
\big(1-p+p\Fbar( \tau_{F}-v/b_n)\big)  \Gbar(\tau_{F}-v/b_n)
- (1-p) \int_{\tau_{F}-v/b_n}^{\tau_{G}-u/a_n}  \rmd G(z)\cr
&& \hskip6cm 
- p \int_{\tau_{F}-v/b_n}^{\tau_{F}}  \Fbar (z)  \rmd G(z).
\eea
The last equality comes about because $\Fbar(z)=0$ for $z\ge \tau_{F}$.
After a cancellation and multiplying through by $n$,
 the last expression takes the form 
\bea\label{ap6z}
&&
np\Fbar( \tau_{F}-v/b_n) \Gbar(\tau_{F}-v/b_n)
+n(1-p)  \Gbar(\tau_{G}-u/a_n)     \cr
&&\hskip4cm 
-n p \int_{\tau_{F}-v/b_n}^{\tau_{F}}  \Fbar (z)  \rmd G(z).
\eea
After integrating by parts and use of the mean value theorem, \eqref{ap6z} becomes 
\bea\label{ap66}
&&
n(1-p)  \Gbar(\tau_{G}-u/a_n)
+n p \int_{\tau_{F}-v/b_n}^{\tau_{F}}  \Gbar (z)  \rmd F(z)\cr
&&=
n(1-p)  \Gbar(\tau_{G}-u/a_n)
+ np  \Gbar(\tau_{F}-z_n) \Fbar(\tau_{F}-v/b_n),
\eea
where $0 \le z_n \le v/b_n$.
For the same   $a_n$ and $b_n$ 
 the expression on the RHS of \eqref{ap66} has limit
\be\label{ap9}
(1-p) u^\gamma +   p \Gbar(\tau_{F})v^\beta .
\ee
Putting together \eqref{ap4}, \eqref{ap6} and \eqref{ap9} and recalling \eqref{ap1} and \eqref{ap2} gives 
\bean 
&&
\lim_{n\to\infty} 
P\big( a_n(\tau_H- M(n))\le u, b_n(\tau_J- M_u(n))\le v \big) \cr
&&\cr
&&= 
1+e^{-p \Gbar(\tau_{F})v^\beta - (1-p)u^\gamma} 
 - e^{-p \Gbar(\tau_{F})v^\beta }
- e^{- (1-p)u^\gamma}, \ \ 
\eean
and hence the  independent distributions in \eqref{as1}.

\medskip
\noindent {\bf Case 2:}\  
 Assume  $\tau_{F}<\tau_G<\infty$ and $p=1$,
so that $F^*\equiv F$, $\Fbar^*=\Fbar$  and  $\tau_J=\tau_{F}= \tau_{F^*}=\tau_H<\tau_G<\infty$.
Assume also \eqref{F^*GD2}.
In this case we choose just one norming sequence $a_n$ so that 
 $n\Gbar(\tau_G-1/a_n)\to 1$, thus
 $a_n\sim (naL(1/a_n)) ^{1/\beta} $, and the component parts on the RHS of \eqref{ap1} can be calculated as
\bea\label{ap14}  
&&
n\Hbar( \tau_{F}-u/a_n) 
=
n\Fbar( \tau_{F}-u/a_n)\Gbar( \tau_{F}-u/a_n)\cr
&&\cr
&&=
 n a (1+o(1)) \frac{u^\beta L(u/a_n)}{a_n^\beta}
 (\Gbar( \tau_{F})+o(1))
\to \Gbar( \tau_{F}) u^\beta.
\quad
\eea
The asymptotic relations follow from \eqref{F^*GD2}.
Similarly, using \eqref{ap5}, 
\be\label{ap15}  
n\Jbar( \tau_{F}-v/a_n) 
\to \Gbar( \tau_{F})v^\beta.
\ee

To deal with the second summand in \eqref{ap2} in this case, we have to consider two subcases. We set $b_n=a_n$ for this part.
Suppose first that $0<u<v$.
Working from the second line of \eqref{ap5a}, we look at
\be\label{ap17}
 \Hbar( \tau_{F}-v/b_n)
-\int_{\tau_{F}-v/b_n}^{\tau_{F}}   \Fbar (z)  \rmd G(z)
=
\int_{\tau_{F}-v/a_n}^{\tau_{F}}   \Gbar(z)  \rmd F(z).
\ee
This follows from an integration by parts.
Using the mean value theorem, the RHS here is 
\be\label{ap17a}
 \Gbar(\tau_{F}-z_n) \Fbar(\tau_{F}-v/a_n)
=
a (1+o(1)) \frac{v^\beta L(v/a_n)}{a_n^\beta}
\big( \Gbar(\tau_{F}) +o(1)\big)
\ee
where $0 \le z_n \le v/a_n$.
Multiplying \eqref{ap17a}  by $n$ and letting $n\to\infty$ we obtain the limit
$\Gbar(\tau_{F}) v^\beta$ on the RHS, 
and hence
\be\label{ap18}
\lim_{n\to\infty} 
\Big(\int_{\tau_J-v/a_n}^{\tau_H-u/a_n}  \Fbar^* (z)\rmd G(z)+
H( \tau_{F}-v/a_n)
\Big)^n 
=e^{-v^\beta\Gbar(\tau_{F})}.
\ee
Thus, via \eqref{ap1}, putting together \eqref{ap14}, \eqref{ap15} and \eqref{ap18} we get 
\be\label{ap18a}
\lim_{n\to\infty}
P\big( a_n(\tau_H- M(n))\le u,\  a_n(\tau_J- M_u(n))\le v \big) 
=1- e^{-\Gbar(\tau_{F})u^\beta}.
\ee
The second  subcase is when $u>v$, so
$\tau_H-u/a_n <\tau_J-v/a_n$. In this case we use the 
second line on the RHS of \eqref{JD} to write
\bean
&&
P\big(M(n)< \tau_H-u/a_n,\  M_u(n)<\tau_J-v/a_n \big)
=
H^n(  \tau_{F}-u/a_n)\cr
&&=
 \big(1-\Fbar( \tau_{F}-u/a_n) )\Gbar(  \tau_{F}-u/a_n) )\big)^n
\to e^{- \Gbar( \tau_{F})u^\beta},
\eean
by \eqref{ap14}.  So  we get  \eqref{ap18a} again and hence \eqref{as31}.

\medskip
\noindent {\bf Case 3:}\  Assume $\tau_G<\tau_{F}$ and  keep $0<p<1$.
In this case $\tau_J=\tau_H=\tau_G< \tau_{F} <\tau_{F^*}=\infty$.
Assume also the first relation in \eqref{F^*GD}.
Once again we calculate the component parts on the RHS of \eqref{ap1}.
This time we have to scale differently.
We take $a_n$ to satisfy  $n\Gbar(\tau_G-1/a_n)\to 1$ and 
 $a_n\sim (na_GL_G(1/a_n)) ^{1/\gamma} $
as in Case 1, but set
$b_n:= \sup\{x>1/\tau_G: x^{-1}\Gbar(\tau_G-1/x)\ge 1/n\}$.
Then $n\Gbar(\tau_G-1/b_n)\sim b_n\to\infty$ and 
$b_n\sim (na_GL_{G}(1/b_n))^{\frac{1}{1+\gamma}}$.
Thus $a_n$ is bigger than $b_n$, ultimately, and, for any $u,v>0$,
 we have $u/a_n<v/b_n$  hence $\tau_G-u/a_n>\tau_G-v/b_n$ for $n$ large enough.

In place of \eqref{ap3} we have,
by continuity of $F$, the first relation in  \eqref{F^*GD},
and the relation  $a_n^\gamma\sim na_GL_G(1/a_n)$,
\be\label{ap10}  
n\Hbar( \tau_H-u/a_n) 
=
n\big(1-p F( \tau_G-u/a_n)\big)\Gbar( \tau_G-u/a_n)
\to \big(1-p F( \tau_G))u^\gamma,
\ee
For the $J$ part, modify \eqref{ap5} to 
\bea\label{ap11a}  
&&
\Jbar( \tau_J-v/b_n) 
=p \int_{\tau_G-v/b_n}^{\tau_G}  \Gbar (z)\rmd F(z)=
 p\int_0^{v/b_n}  \Gbar(\tau_G-z)  f(\tau_G-z)  \rmd z    \cr
&&\cr
&&=
 p\int_0^{v/b_n}(1+o(1))  a_Gz^\gamma L_G(z) f(\tau_G-z)  \rmd z   \cr
&&=
p  a_G
 f(\tau_G-z_n)  \int_0^{v/b_n}   (1+o(1)) z^\gamma L_G(z) \rmd z,
\eea 
where $0 \le z_n\le v/b_n$. To deal with $\Gbar$ we used  the first relation in \eqref{F^*GD} and to deal with $F$ we used  the mean value  theorem, together with the assumption that,  in a neighbourhood of $\tau_G$, $F$  has a density $f$ which is positive and continuous at $\tau_G$.
 Letting $n\to \infty$ with 
 $b_n ^{1+\gamma}   \sim na_GL_{G}(1/b_n)$, we get from \eqref{ap11a} and  the slow variation of $L_G$ that
\be\label{ap11b}
\lim_{n\to\infty} n \Jbar( \tau_J-v/b_n) 
=
\frac{p f(\tau_G) v^{1+\gamma}}{1+\gamma}.
\ee

Finally, in place of \eqref{ap5a}: 
 for any $u,v>0$ we have  $\tau_G-u/a_n>\tau_G-v/b_n$ 
for $n$ large enough.
Use  integration by parts, 
$\Hbar (x)
=\big(1-p+p\Fbar(x)\big)\Gbar(x)$ by \eqref{13},
$\tau_H=\tau_J=\tau_G$,  and the mean value theorem, to calculate
\bean 
&&
n\Hbar( \tau_J-v/b_n)-n  \int_{\tau_J-v/b_n}^{\tau_H-u/a_n}\big(1-p  F (z)\big)\rmd G(z)\cr
&&\cr
&&=
n\big(1-pF( \tau_{G}-u/a_n)\big)    \Gbar(\tau_G-u/a_n)
+np  \int_{\tau_G-v/b_n}^{\tau_G-u/a_n}  \Gbar(z) \rmd F(z) \cr
&&=
\big(1-pF( \tau_{G})+o(1)\big) (1+o(1)) u^\gamma\cr
&&  \hskip4cm 
+n p f(\tau_G-z_n)   \int_{u/a_n}^{v/b_n}  (1+o(1)) a_Gz^\gamma L_G(z)\rmd  z \cr
&&
=
\big(1-pF( \tau_{G})\big)   u^\gamma
+\frac{p ( f(\tau_G)+o(1) )}{1+\gamma}
\Big(v^{1+\gamma} -  O\Big(\frac{b_n^{1+\gamma}u^{1+\gamma}  }{ a_n^{1+\gamma}} \Big)\Big),
\eean
where $u/a_n \le z_n \le v/b_n$.
With the previous choices of $a_n$ and  $b_n $ we find from this that 
\bea\label{ap13}
&&
\lim_{n\to\infty}
n
\Big(\Hbar( \tau_J-v/b_n)-  \int_{\tau_J-v/b_n}^{\tau_G-u/b_n}  \Fbar^* (z)\rmd G(z)\Big)\cr
&&
=
\big(1-pF( \tau_{G})\big)u^\gamma
+  \frac{p f(\tau_G) v^{1+\gamma}}{1+\gamma} .
\eea
Put together \eqref{ap10}, \eqref{ap11b} and \eqref{ap13} and recall \eqref{ap1} and \eqref{ap2} to get \eqref{as2}.    

For this case, we can easily check that all the working  
so far remains true under the same assumptions when $p=1$, and for none of it  
is the value of $F$ beyond $\tau_G$ relevant, so the results
hold equally well  when $\tau_{F}=\infty$.
\halmos


\bibliography{bibfile.bib}


\newpage
\section{Illustrations of cdfs}\label{ssfigs}
In this section we plot the cdfs of some of the distributions in Section \ref{s1}.
Most figures have  $F=U[0,a]$ and $G=U[0,10]$, with  $a=15,10,5$.
In this scenario,   $a=15$ represents insufficient follow-up, $a=10$ is a marginal case (insufficient), and $a=5$ represents sufficient follow-up.
All plots were done using the statistical package R.
%
%


%

\begin{figure}[h]
\centering
\begin{minipage}{.3\linewidth}
\centering
\includegraphics[width=\linewidth]{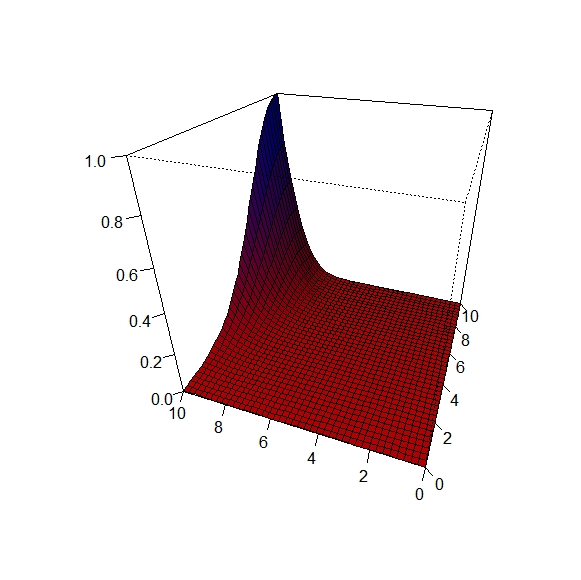}
\end{minipage}\hfill
\begin{minipage}{.3\linewidth}
\centering
\includegraphics[width=\linewidth]{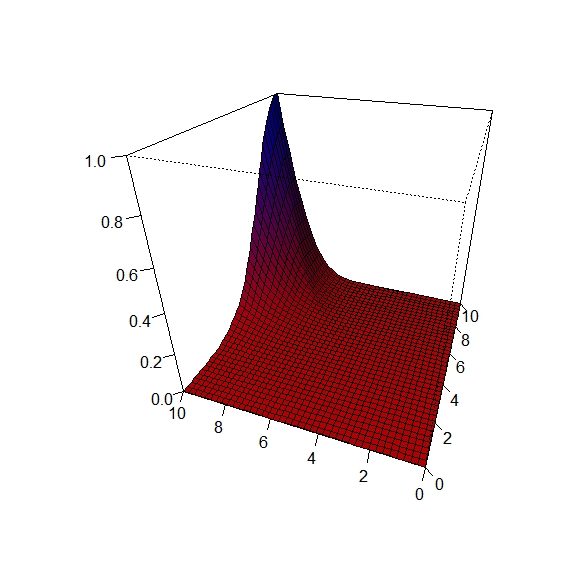}
\end{minipage}\hfill
\begin{minipage}{.3\linewidth}
\centering
\includegraphics[width=\linewidth]{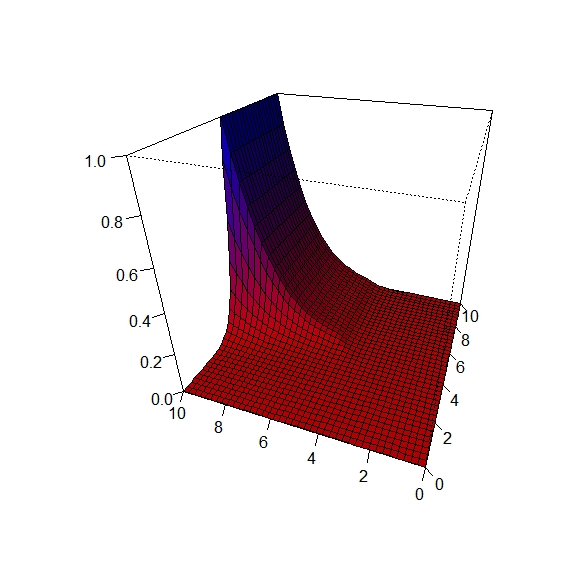}
\end{minipage}\hfill
\caption{The joint cdf of  $M(n)$ and $M_u(n) $ as given by \eqref{JD} with\\ $F=U[0,a]$, $G=U[0,10]$, $p=0.7$. Left to right:  $a=15,10,5$. }\label{fig:1}
\end{figure}

\begin{figure}[h]
\centering
\begin{minipage}{.3\linewidth}
\centering
\includegraphics[width=\linewidth]{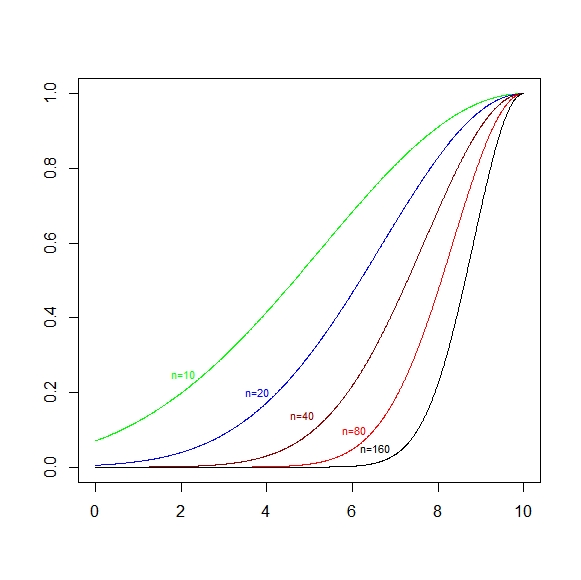}
\end{minipage}\hfill
\begin{minipage}{.3\linewidth}
\centering
\includegraphics[width=\linewidth]{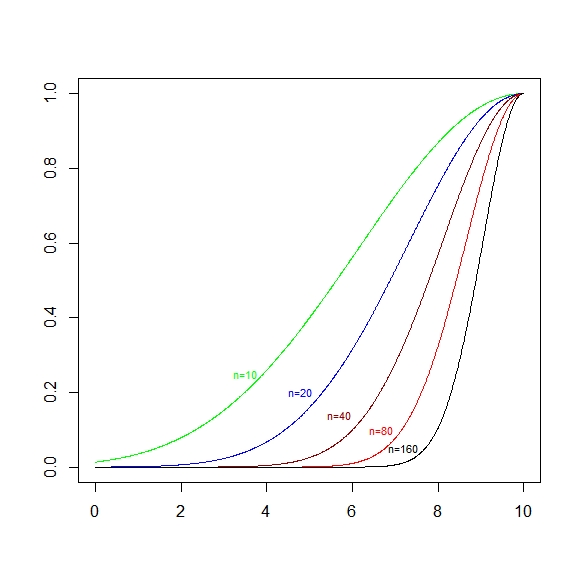}
\end{minipage}\hfill
\begin{minipage}{.3\linewidth}
\centering
\includegraphics[width=\linewidth]{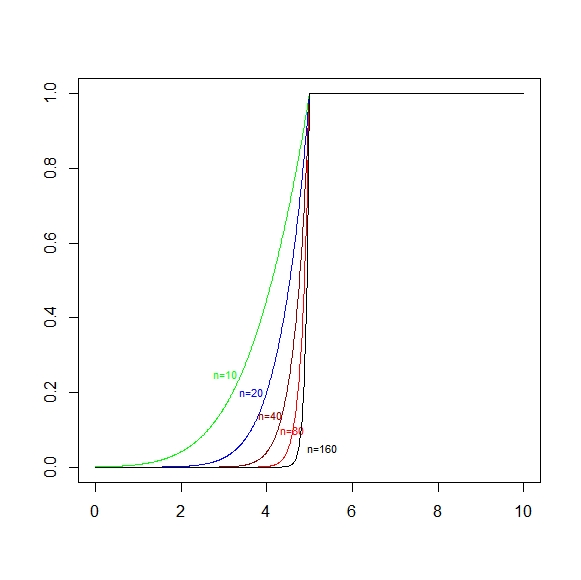}
\end{minipage}\hfill
\caption{The cdf of   $M_u(n) $ as given by \eqref{4f} with  $F=U[0,a]$, $G=U[0,10]$, $p=0.7$. Left to right: $a=15, 10, 5$.}\label{fig:2}
\end{figure}

\begin{figure}[h]
\centering
\begin{minipage}{.3\linewidth}
\centering
\includegraphics[width=\linewidth]{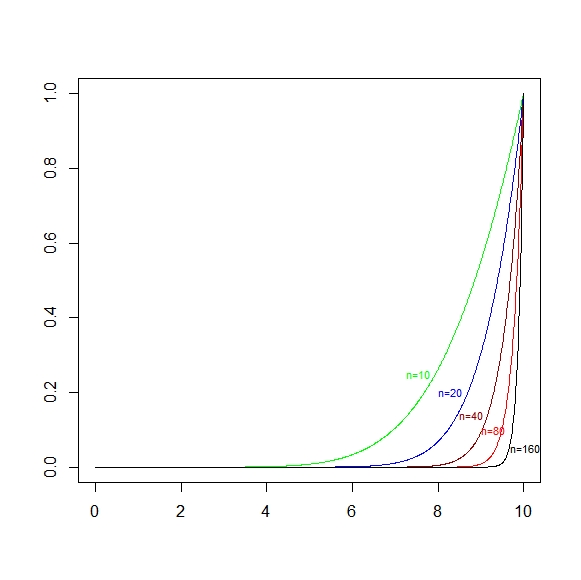}
\end{minipage}\hfill
\begin{minipage}{.3\linewidth}
\centering
\includegraphics[width=\linewidth]{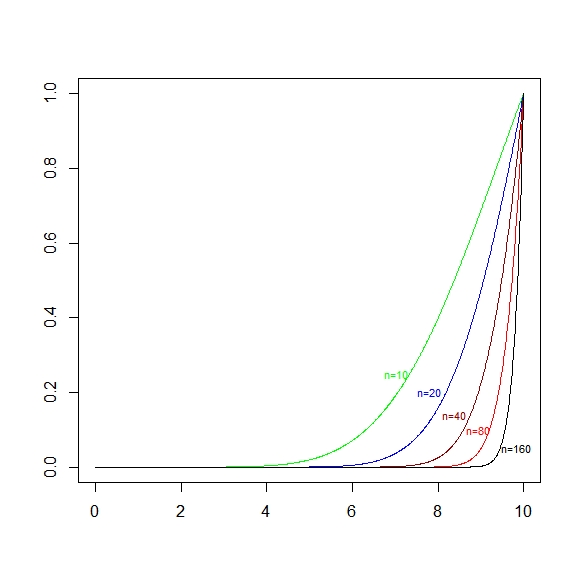}
\end{minipage}\hfill
\begin{minipage}{.3\linewidth}
\centering
\includegraphics[width=\linewidth]{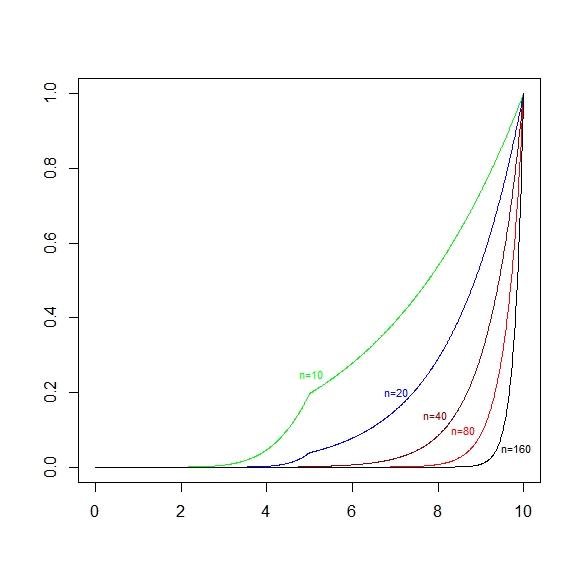}
\end{minipage}\hfill
\caption{The cdf of   $M(n)$ as given by \eqref{4h} with  $F=U[0,a]$, $G=U[0,10]$, $p=0.7$. Left to right: $a=15,10,5$.}\label{fig:3}
\end{figure}

\begin{figure}[h]
\centering
\begin{minipage}{.3\linewidth}
\centering
\includegraphics[width=\linewidth]{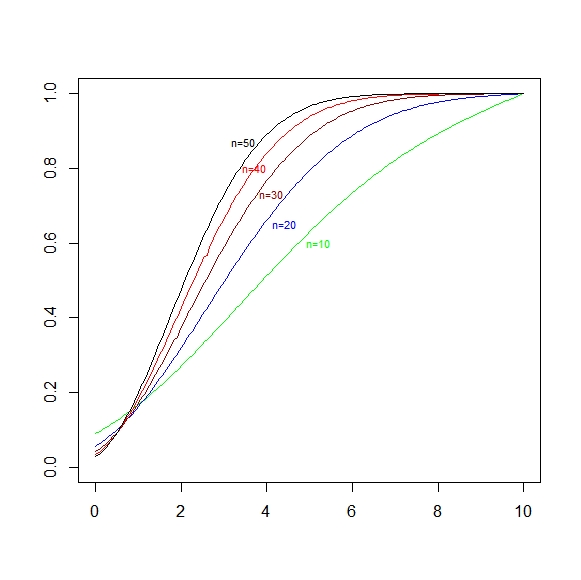}
\end{minipage}\hfill
\begin{minipage}{.3\linewidth}
\centering
\includegraphics[width=\linewidth]{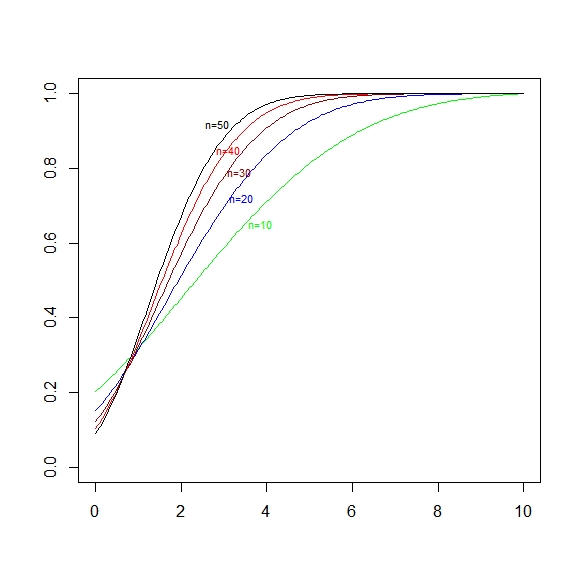}
\end{minipage}\hfill
\begin{minipage}{.3\linewidth}
\centering
\includegraphics[width=\linewidth]{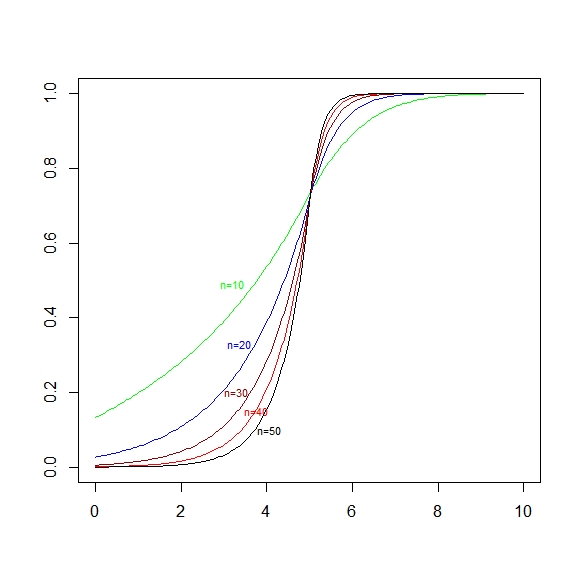}
\end{minipage}\hfill
\caption{The cdf of   $M(n)-M_u(n) $ as given by \eqref{12a} with
 $F=$U[0, a$]$, $G=U[0,10]$,  $p=0.7$. Left to right: $a=15,10,5$.}\label{fig:4}
\end{figure}

\begin{figure}[ht]
\centering
\begin{minipage}{.3\linewidth}
\centering
\includegraphics[width=\linewidth]{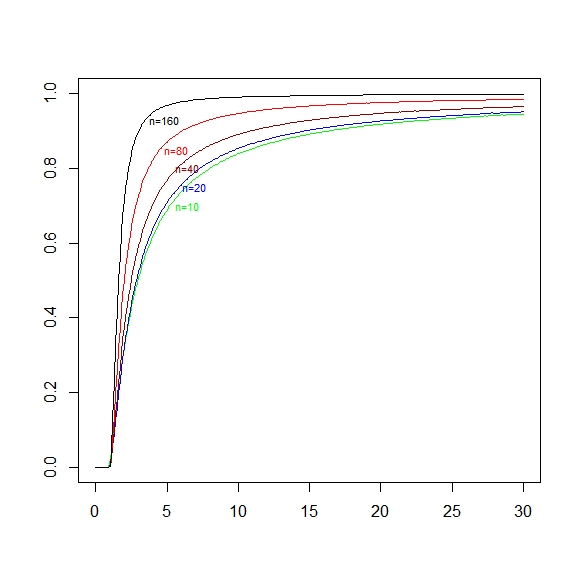}
\end{minipage}\hfill
\begin{minipage}{.3\linewidth}
\centering
\includegraphics[width=\linewidth]{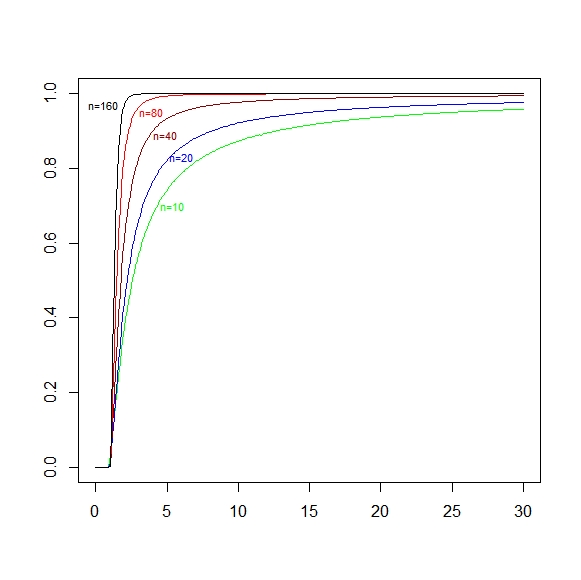}
\end{minipage}\hfill
\begin{minipage}{.3\linewidth}
\centering
\includegraphics[width=\linewidth]{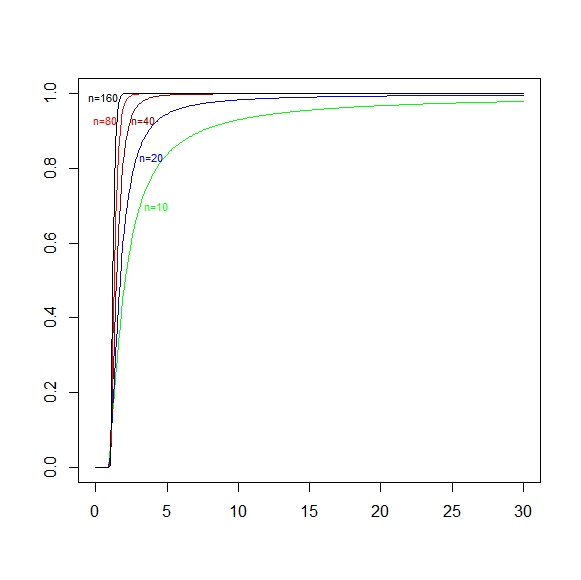}
\end{minipage}\hfill
\caption{The cdf of   $M(n)/M_u(n) $ (conditional on $M_u(n) >0$)  as given by \eqref{12b} with 
 $F=$Exp($\lambda=0.0802$), $G=U[0,10]$, $p=0.7$.}\label{fig:5}
\end{figure}


\begin{figure}[h]
\centering
\begin{minipage}{.3\linewidth}
\centering
\includegraphics[width=\linewidth]{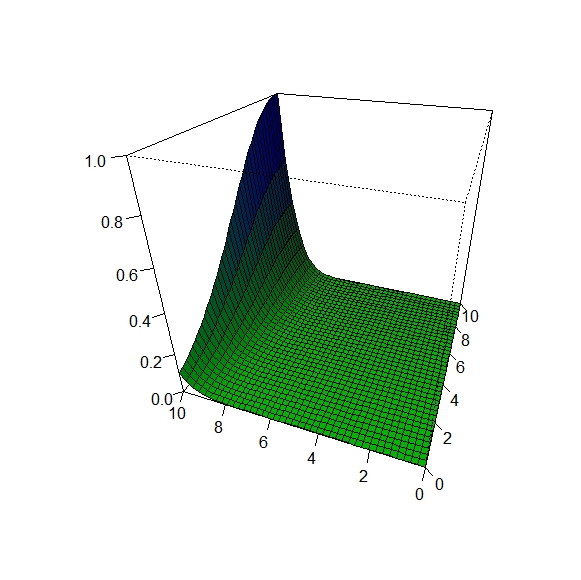}
\vspace{1ex}
\end{minipage}\hfill
\begin{minipage}{.3\linewidth}
\centering
\includegraphics[width=\linewidth]{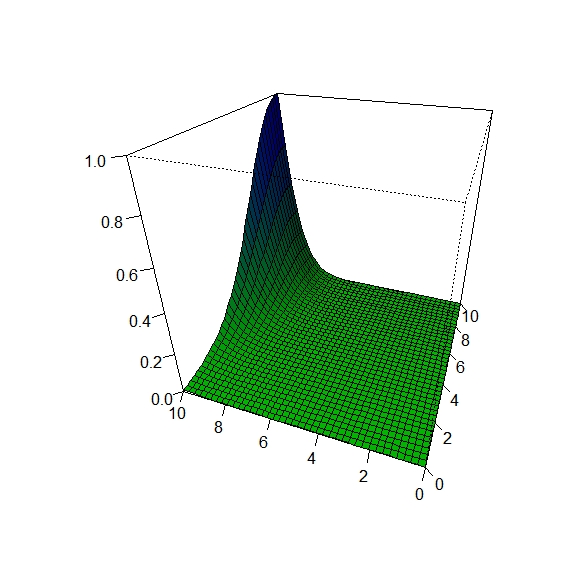}
\vspace{1ex}
\end{minipage}\hfill
\begin{minipage}{.3\linewidth}
\centering
\includegraphics[width=\linewidth]{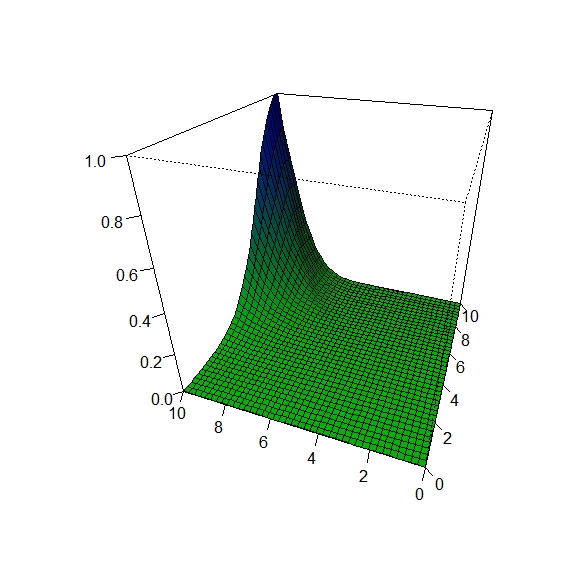}
\vspace{1ex}
\end{minipage}\hfill
\centering
\begin{minipage}{.3\textwidth}
\centering
\includegraphics[width=\linewidth]{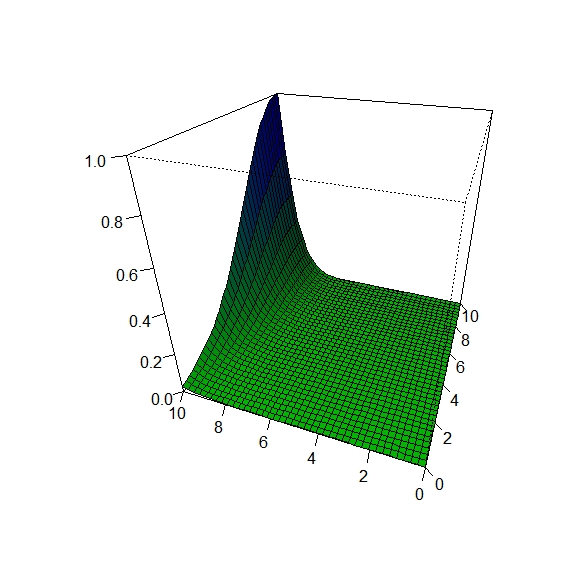}
\vspace{1ex}
\end{minipage}\hfill
\begin{minipage}{.3\linewidth}
\centering
\includegraphics[width=\linewidth]{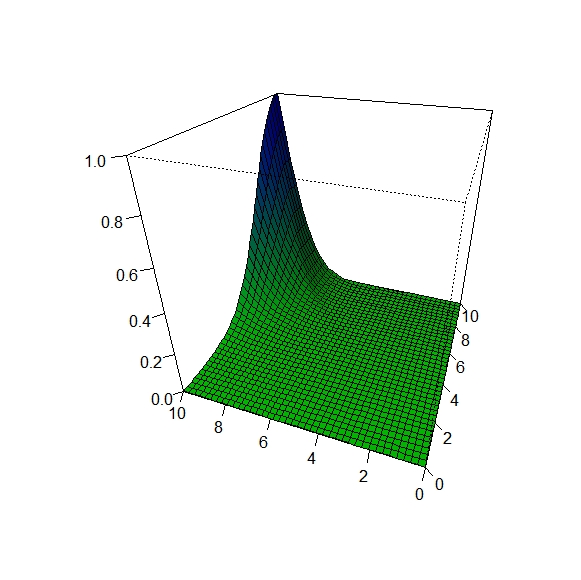}
\vspace{1ex}
\end{minipage}\hfill
\begin{minipage}{.3\linewidth}
\centering
\includegraphics[width=\linewidth]{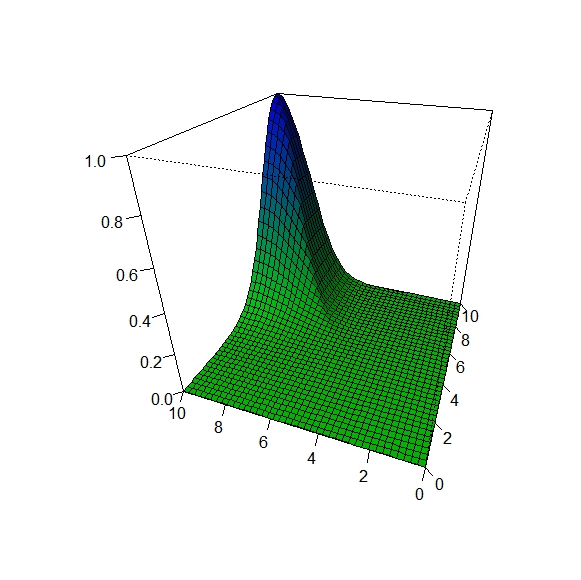}
\vspace{1ex}
\end{minipage}\hfill
\begin{minipage}{.3\linewidth}
\centering
\includegraphics[width=\linewidth]{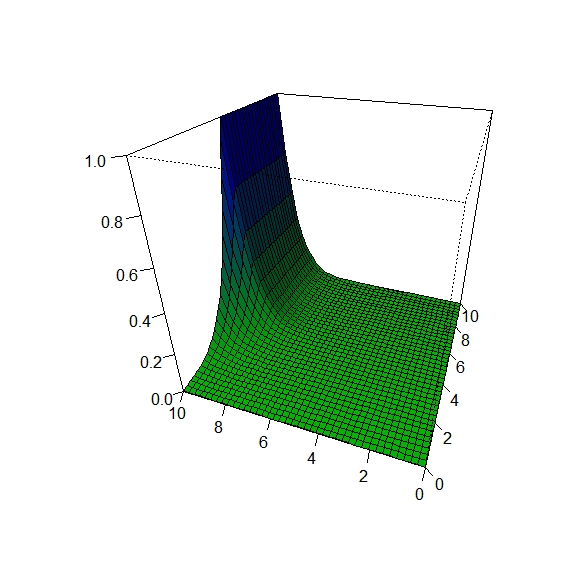}
\vspace{1ex}
\end{minipage}\hfill
\begin{minipage}{.3\linewidth}
\centering
\includegraphics[width=\linewidth]{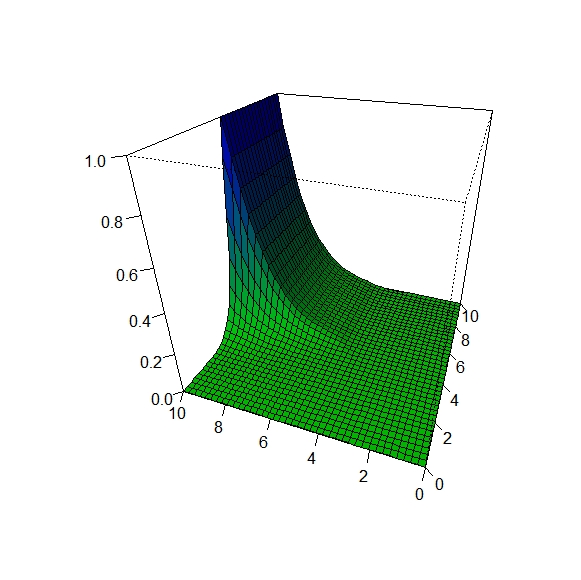}
\vspace{1ex}
\end{minipage}\hfill
\begin{minipage}{.3\linewidth}
\centering
\includegraphics[width=\linewidth]{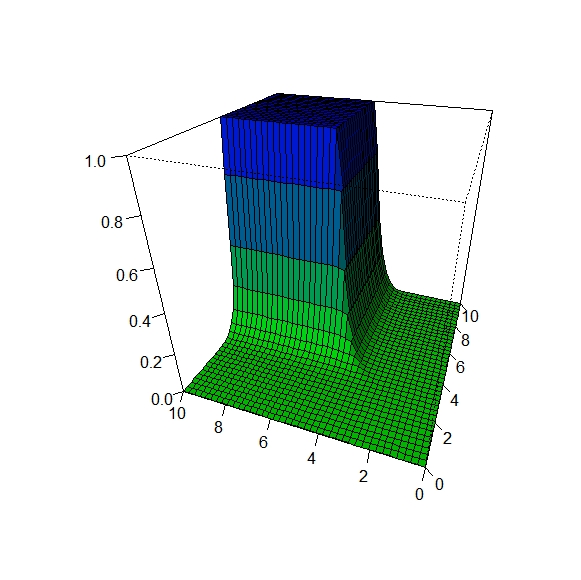}
\vspace{1ex}
\end{minipage}\hfill
\caption{The joint cdf of  $M(n)$ and $M_u(n) $ as given by \eqref{JD} with
 $F=U[0,a]$, $G=U[0,10]$, $n=20$. 
Left to right:  $p=0.33,0.66,1$. 
 Top to bottom: $a=15,10,5$.}\label{fig:6}
\end{figure}

%

\begin{figure}[h]
\centering
\begin{minipage}{.3\linewidth}
\centering
\includegraphics[width=\linewidth]{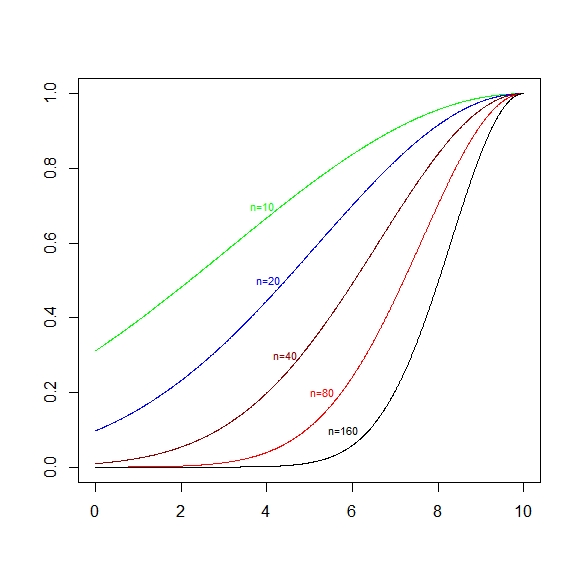}
\vspace{1ex}
\end{minipage}\hfill
\begin{minipage}{.3\linewidth}
\centering
\includegraphics[width=\linewidth]{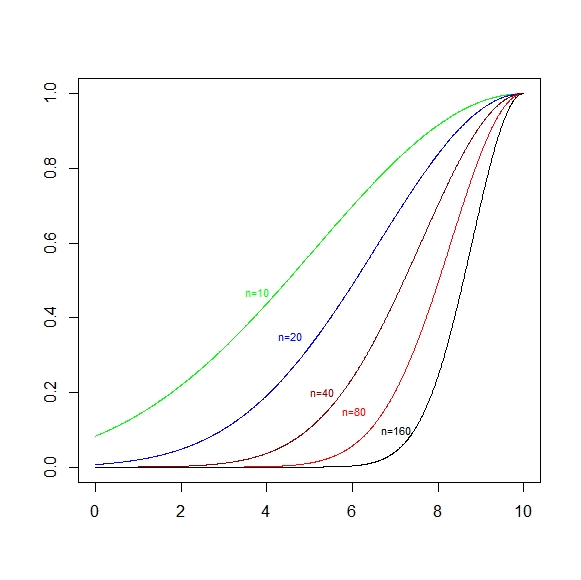}
\vspace{1ex}
\end{minipage}\hfill
\begin{minipage}{.3\linewidth}
\centering
\includegraphics[width=\linewidth]{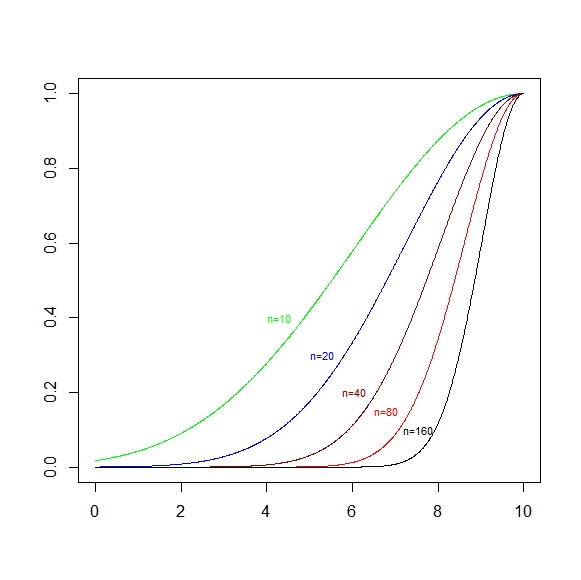}
\vspace{1ex}
\end{minipage}\hfill
\begin{minipage}{.3\textwidth}
\centering
\includegraphics[width=\linewidth]{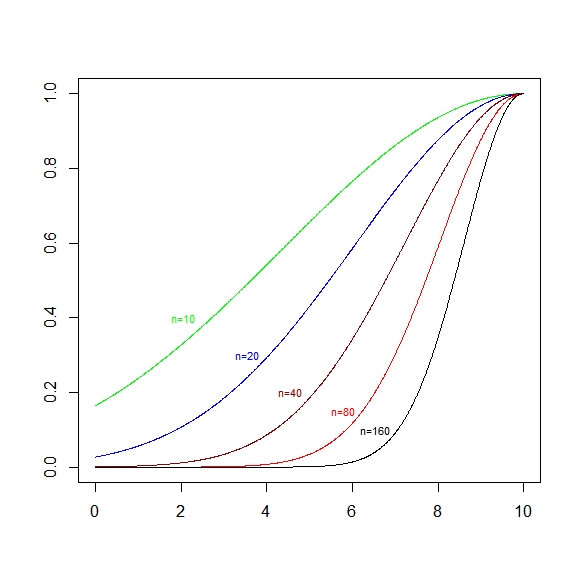}
\vspace{1ex}
\end{minipage}\hfill
\begin{minipage}{.3\linewidth}
\centering
\includegraphics[width=\linewidth]{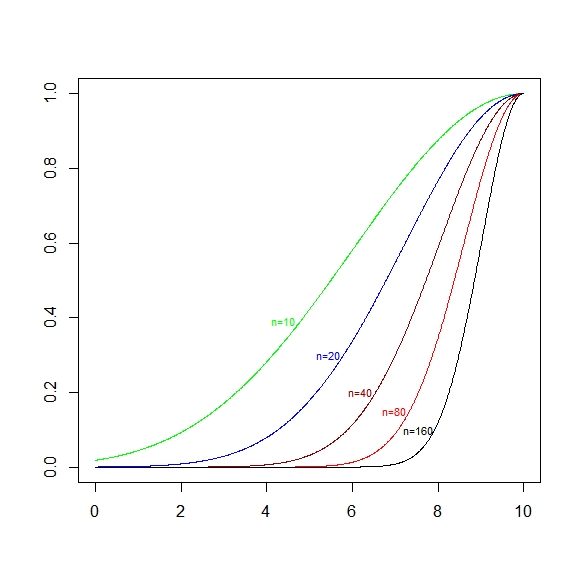}
\vspace{1ex}
\end{minipage}\hfill
\begin{minipage}{.3\linewidth}
\centering
\includegraphics[width=\linewidth]{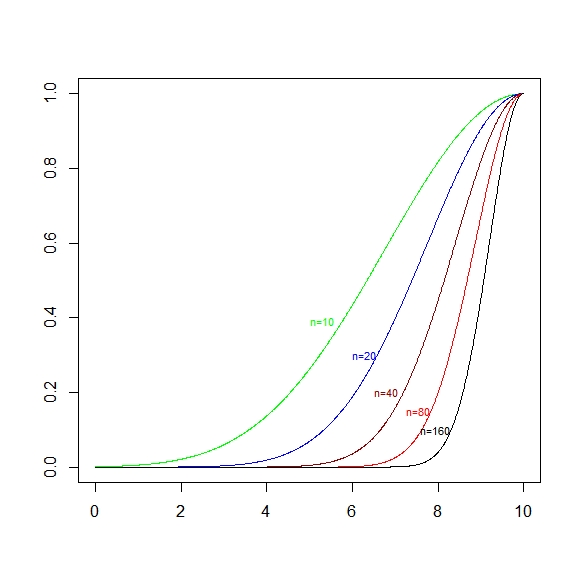}
\vspace{1ex}
\end{minipage}\hfill
\begin{minipage}{.3\linewidth}
\centering
\includegraphics[width=\linewidth]{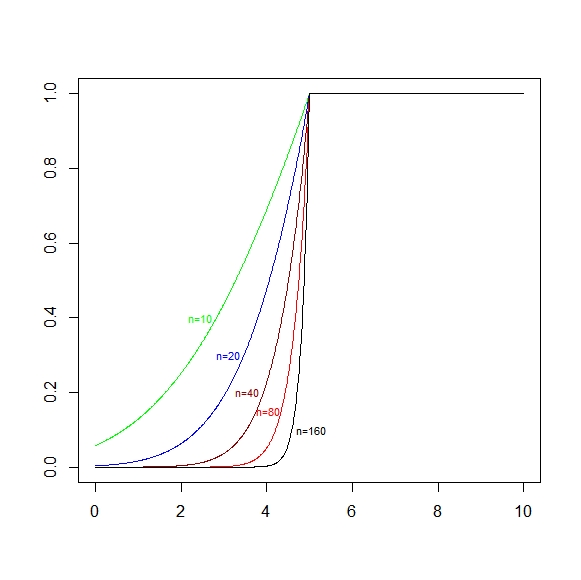}
\vspace{1ex}
\end{minipage}\hfill
\begin{minipage}{.3\linewidth}
\centering
\includegraphics[width=\linewidth]{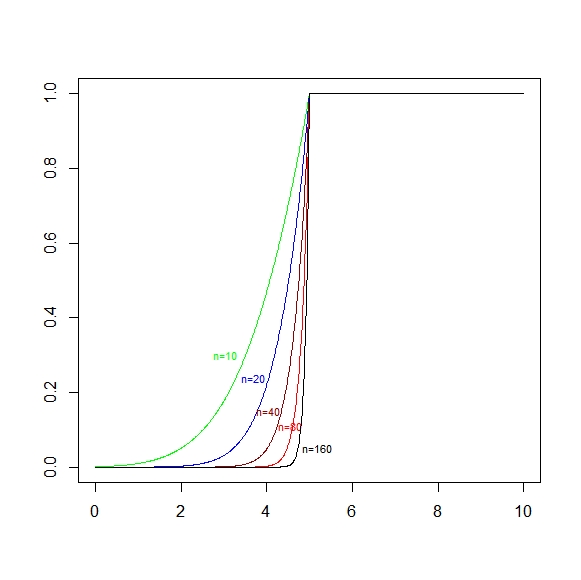}
\vspace{1ex}
\end{minipage}\hfill
\begin{minipage}{.3\linewidth}
\centering
\includegraphics[width=\linewidth]{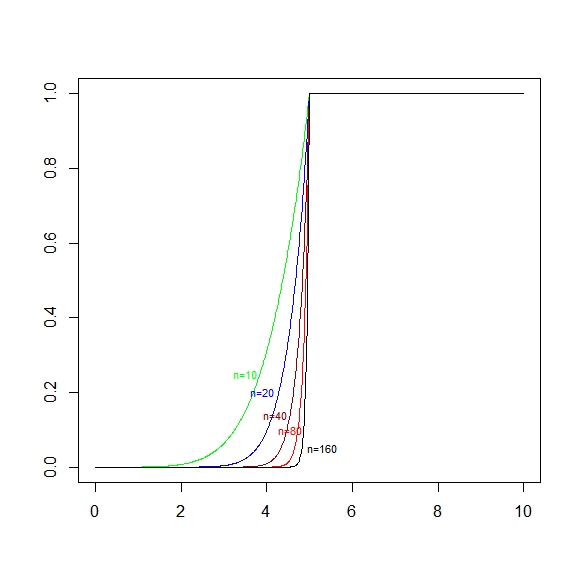}
\vspace{1ex}
\end{minipage}\hfill
\caption{The cdf of   $M_u(n) $ as given by \eqref{4f} with $F=U[0,a]$, $G=U[0,10]$. Left to right: $p=0.33,0.66,1$.  Top to bottom: $a=15,10,5$. }\label{fig:7}
\end{figure}


%

\begin{figure}[h]
\centering
\begin{minipage}{.3\linewidth}
\centering
\includegraphics[width=\linewidth]{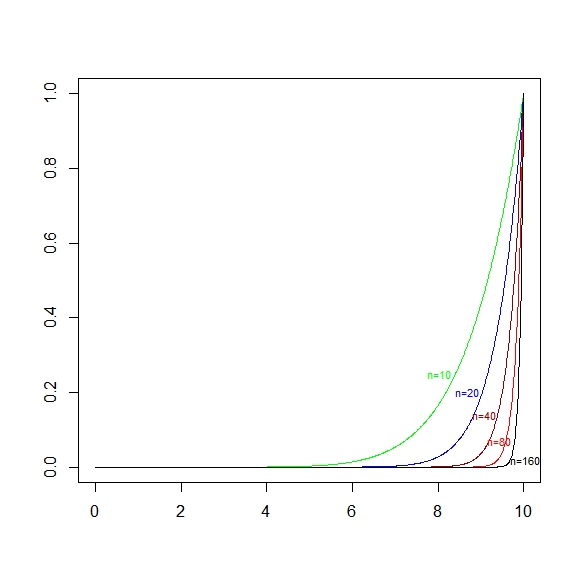}
\vspace{1ex}
\end{minipage}\hfill
\begin{minipage}{.3\linewidth}
\centering
\includegraphics[width=\linewidth]{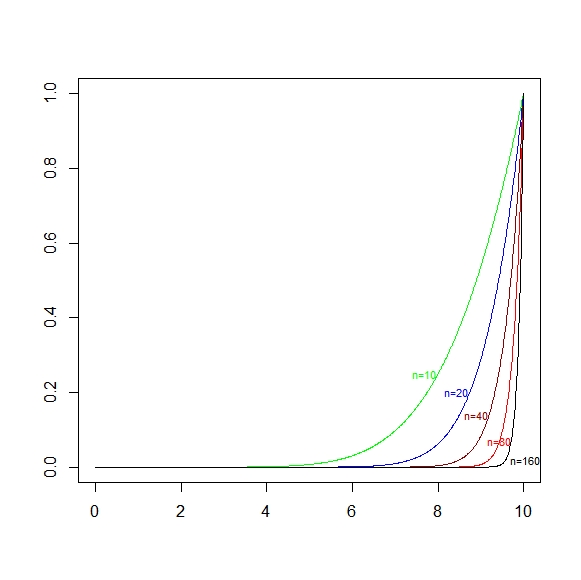}
\vspace{1ex}
\end{minipage}\hfill
\begin{minipage}{.3\linewidth}
\centering
\includegraphics[width=\linewidth]{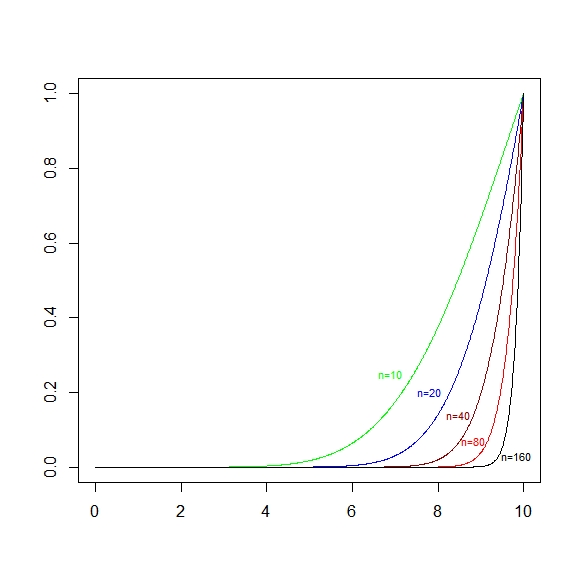}
\vspace{1ex}
\end{minipage}\hfill
\begin{minipage}{.3\textwidth}
\centering
\includegraphics[width=\linewidth]{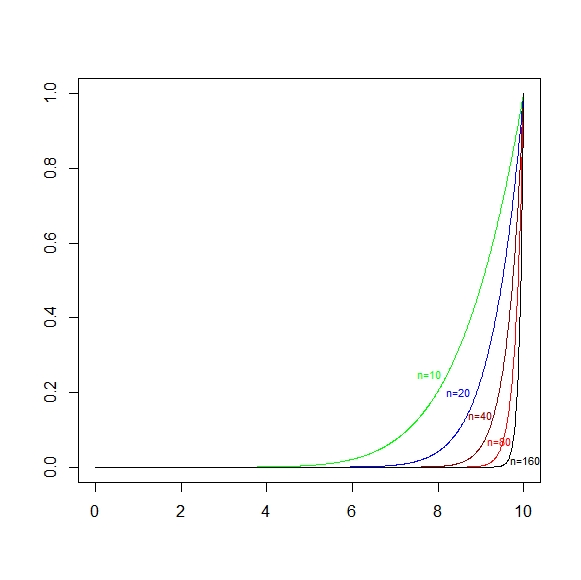}
\vspace{1ex}
\end{minipage}\hfill
\begin{minipage}{.3\linewidth}
\centering
\includegraphics[width=\linewidth]{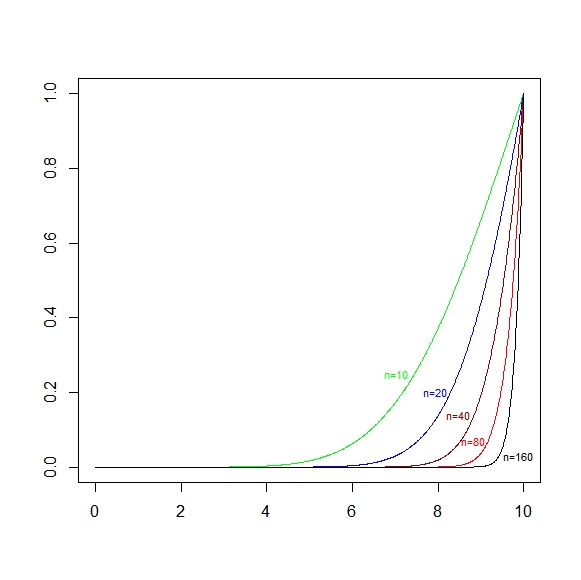}
\vspace{1ex}
\end{minipage}\hfill
\begin{minipage}{.3\linewidth}
\centering
\includegraphics[width=\linewidth]{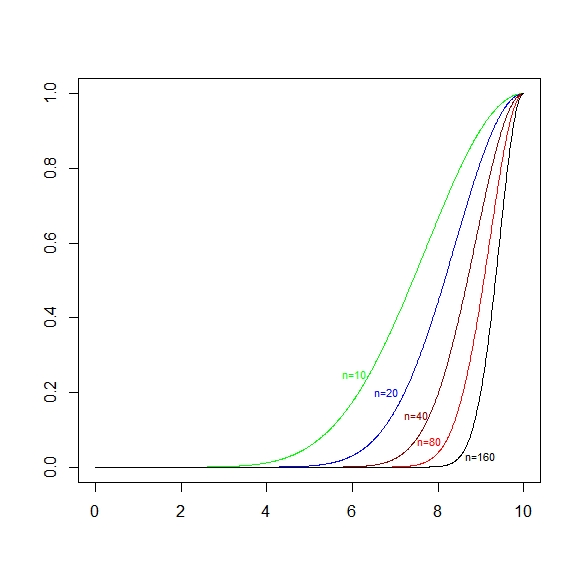}
\vspace{1ex}
\end{minipage}\hfill
\begin{minipage}{.3\linewidth}
\centering
\includegraphics[width=\linewidth]{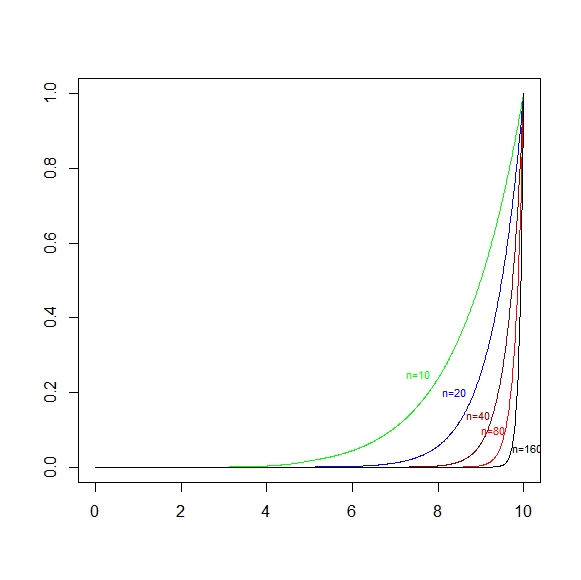}
\vspace{1ex}
\end{minipage}\hfill
\begin{minipage}{.3\linewidth}
\centering
\includegraphics[width=\linewidth]{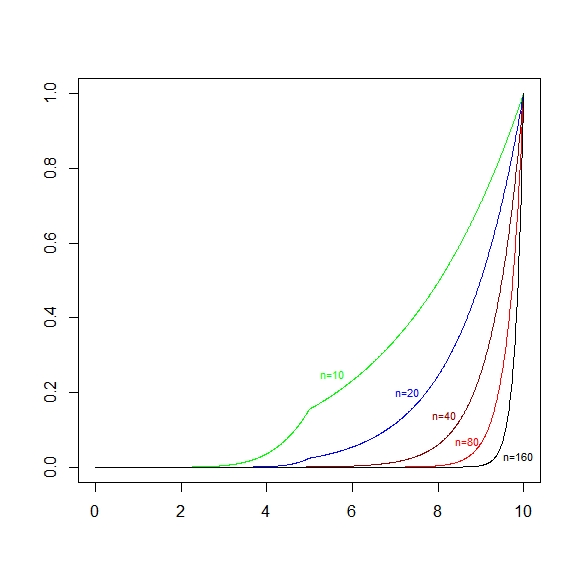}
\vspace{1ex}
\end{minipage}\hfill
\begin{minipage}{.3\linewidth}
\centering
\includegraphics[width=\linewidth]{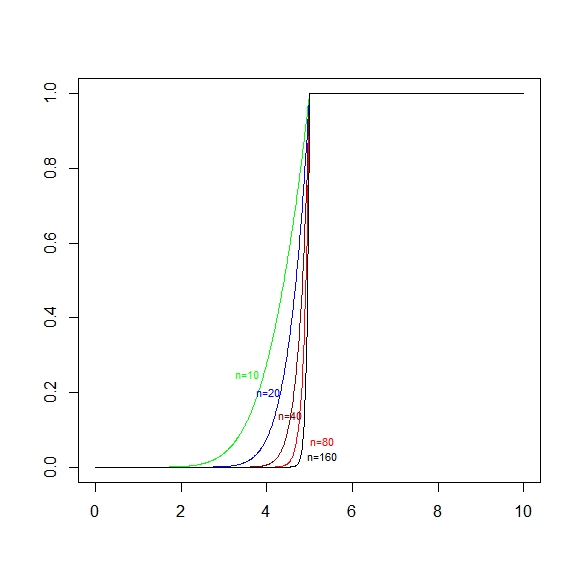}
\vspace{1ex}
\end{minipage}\hfill
\caption{The cdf of   $M(n) $ as given by \eqref{4h} with $F=U[0,a]$, $G=U[0,10]$. Left to right: $p=0.33,0.66,1$.  Top to bottom: $a=15,10,5$. }\label{fig:8}
\end{figure}


%
\begin{figure}[ht]
\centering
\begin{minipage}{.3\linewidth}
\centering
\includegraphics[width=\linewidth]{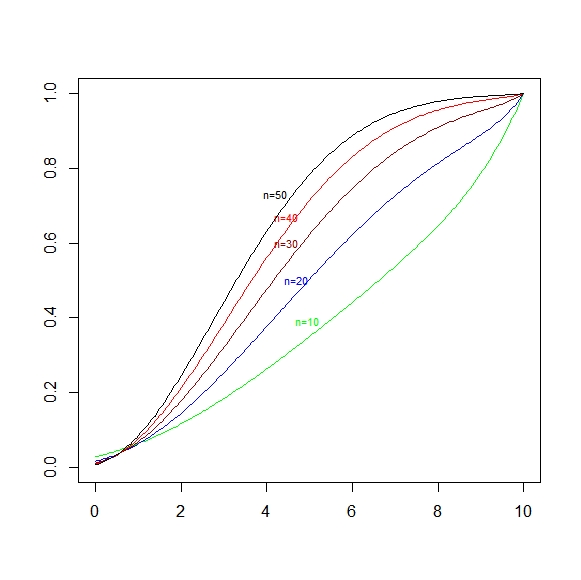}
\vspace{1ex}
\end{minipage}\hfill
\begin{minipage}{.3\linewidth}
\centering
\includegraphics[width=\linewidth]{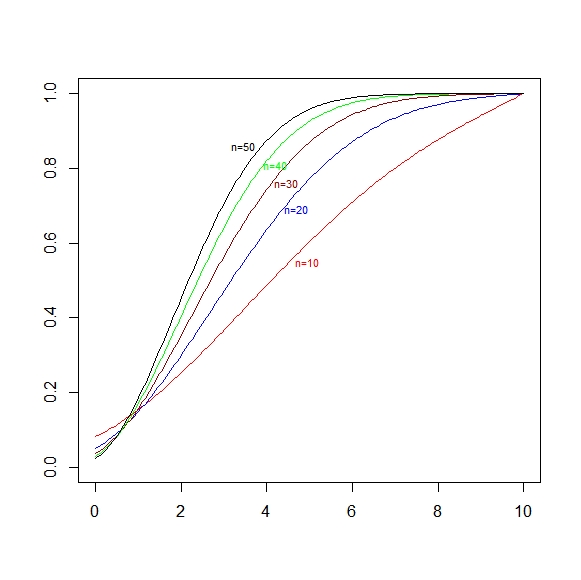}
\vspace{1ex}
\end{minipage}\hfill
\begin{minipage}{.3\linewidth}
\centering
\includegraphics[width=\linewidth]{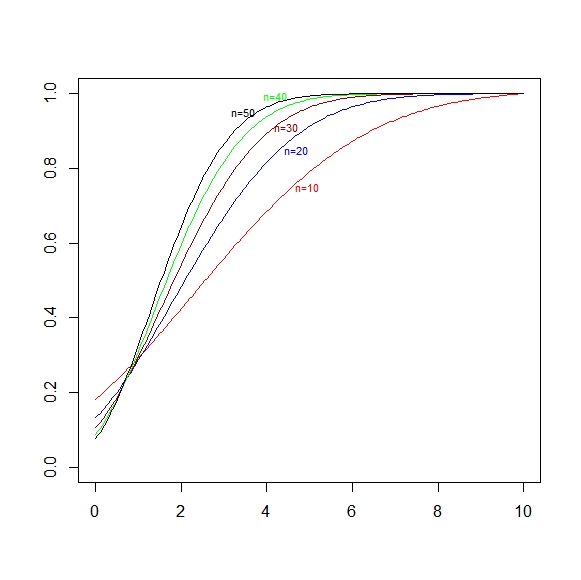}
\vspace{1ex}
\end{minipage}\hfill
\begin{minipage}{.3\textwidth}
\centering
\includegraphics[width=\linewidth]{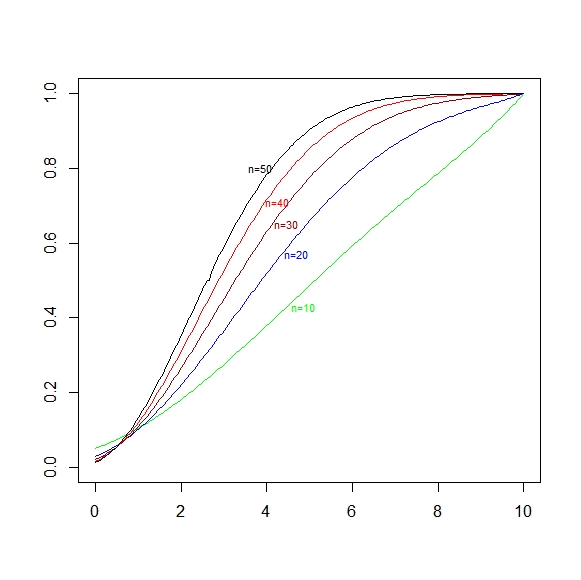}
\vspace{1ex}
\end{minipage}\hfill
\begin{minipage}{.3\linewidth}
\centering
\includegraphics[width=\linewidth]{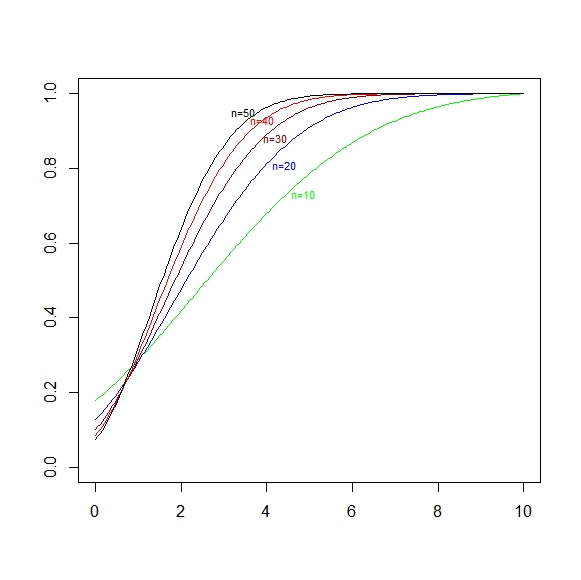}
\vspace{1ex}
\end{minipage}\hfill
\begin{minipage}{.3\linewidth}
\centering
\includegraphics[width=\linewidth]{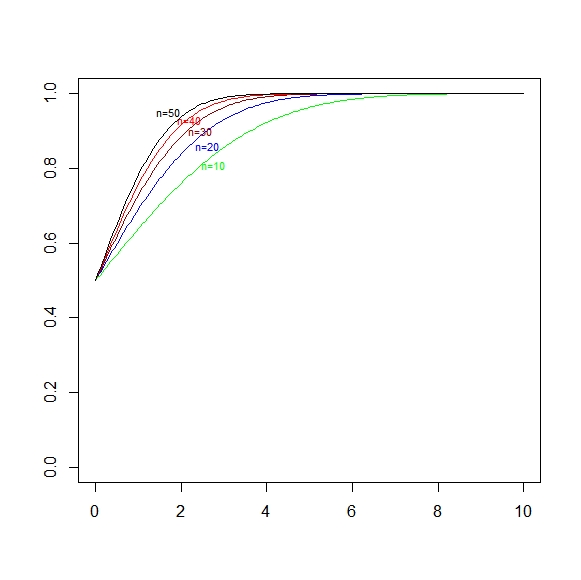}
\vspace{1ex}
\end{minipage}\hfill
\begin{minipage}{.3\linewidth}
\centering
\includegraphics[width=\linewidth]{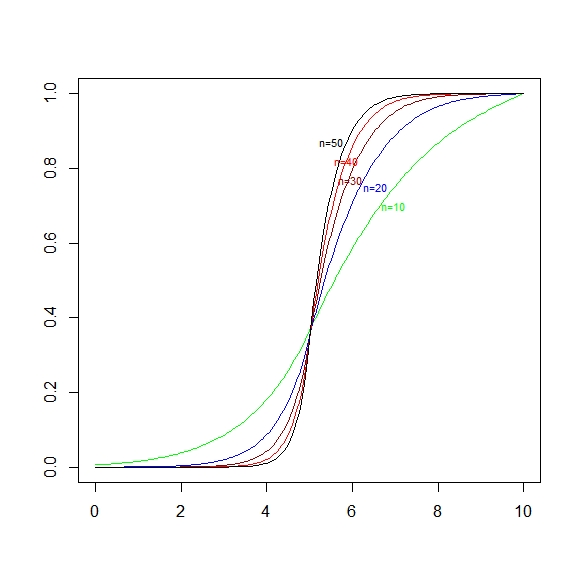}
\vspace{1ex}
\end{minipage}\hfill
\begin{minipage}{.3\linewidth}
\centering
\includegraphics[width=\linewidth]{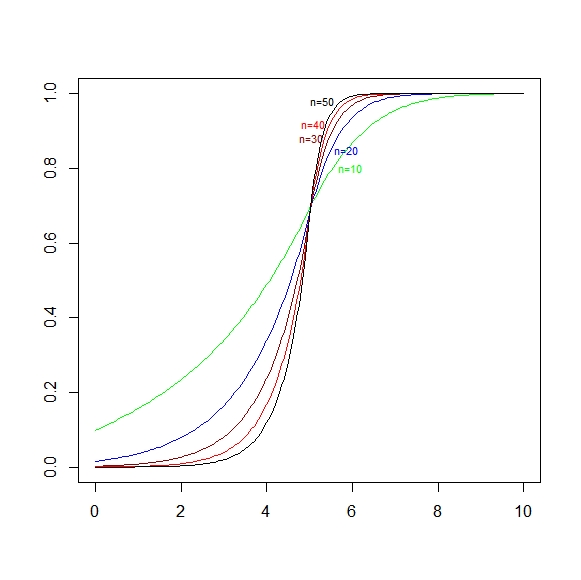}
\vspace{1ex}
\end{minipage}\hfill
\begin{minipage}{.3\linewidth}
\centering
\includegraphics[width=\linewidth]{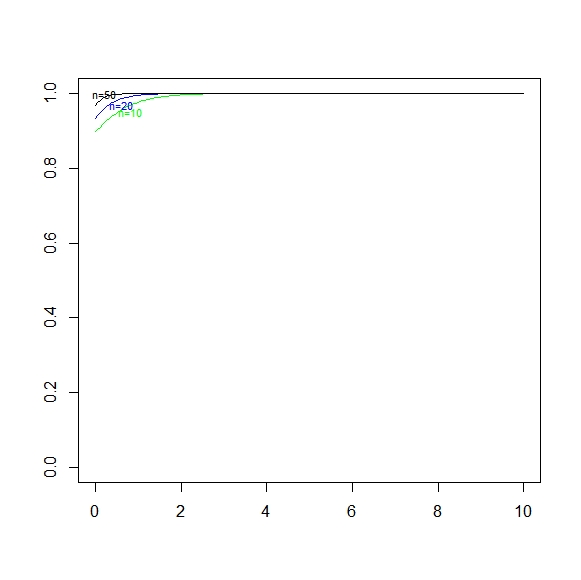}
\vspace{1ex}
\end{minipage}\hfill
\caption{The cdf of   $M(n)-M_u(n) $ as given by \eqref{12a} with $F=U[0,a]$, $G=U[0,10]$.
Left to right: $p=0.33, 0.66, 1$.  Top to bottom: $a=15,10,5$. }\label{fig:9}
\end{figure}


\clearpage

\end{document}